\documentclass[12pt,a4paper]{amsart}
\usepackage{amssymb}
\newtheorem{theorem}[subsection]{Theorem}
\newtheorem{proposition}[subsection]{Proposition}
\newtheorem{corollary}[subsection]{Corollary}
\newtheorem{lemma}[subsection]{Lemma}
\theoremstyle{remark}
\newtheorem*{remark}{Remark}
\numberwithin{equation}{section}
\newcommand{\ffrac}[2]{{\mbox{\large$\frac{#1}{#2}$}}}
\newcommand{\longhookrightarrow}{\lhook\joinrel\relbar\joinrel\rightarrow}
\newcommand{\VB}{V}
\newcommand{\ph}{\varphi}
\newcommand{\rank}{\operatorname{\mathrm{rank}}}
\newcommand{\im}{\operatorname{\mathrm{im}}}
\newcommand{\Hom}{\operatorname{\mathrm{Hom}}}
\newcommand{\End}{\operatorname{\mathrm{End}}}
\newcommand{\cart}{\circledcirc}
\newcommand{\Cart}{\raisebox{-1pt}{\Large$\circledcirc$}}
\newcommand{\on}[2]{\setbox0=\hbox{$#1$}\setbox1=\hbox{$#2$}%
            \dimen0=\wd0\advance\dimen0 by \wd1\divide\dimen0 by 2%
             \ifdim\wd0>\wd1$#1$\hskip-\dimen0$#2$\advance\dimen0 by -\wd1%
              \else$#2$\hskip-\dimen0$#1$\advance\dimen0 by -\wd0%
             \fi%
            \hskip\dimen0}
\newcommand{\onn}[2]{\makebox{\on{#1}{#2}}}
\newcommand{\leftblob}[1]{\hspace*{-1ex}\onn{\raisebox{-.5ex}
     {$\genfrac{}{}{0pt}{}{\hfill\hrulefill}{\hphantom{\hspace*{2ex}#1}}$}}
                {\stackrel{#1}{\bullet}}}
\newcommand{\rightblob}[1]{\onn{\raisebox{-.5ex}
     {$\genfrac{}{}{0pt}{}{\hrulefill\hfill}{\hphantom{\hspace*{2ex}#1}}$}}
                 {\stackrel{#1}{\bullet}}\hspace*{-1ex}}
\newcommand{\middleblob}[1]{\onn{\raisebox{-.5ex}
     {$\genfrac{}{}{0pt}{}{\hrulefill}{\hphantom{\hspace*{2ex}#1}}$}}
                 {\stackrel{#1}{\bullet}}}
\newcommand{\oooo}[4]{\leftblob{#1}\hspace*{-1ex}\middleblob{#2}\hspace*{-1ex}
               \middleblob{#3}\cdots\rightblob{#4}}
\newcommand{\oooospecial}[4]{\leftblob{#1}\hspace*{-1ex}\middleblob{#2}
               \cdots\middleblob{#3}\hspace*{-1ex}\rightblob{#4}}
\newcommand{\ooooo}[5]{\leftblob{#1}\hspace*{-1ex}\middleblob{#2}\hspace*{-1ex}
     \middleblob{#3}\hspace*{-1ex}\middleblob{#4}\cdots\rightblob{#5}}
\newcommand{\son}[4]{\leftblob{#1}\hspace*{-1ex}\middleblob{#2}\hspace*{-1ex}
     \middleblob{#3}\hspace*{-1ex}\middleblob{#4}\cdots}
\begin{document}
\title{Prolongations of Geometric Overdetermined Systems}
\author{Thomas Branson}
\address{Department of Mathematics, University of Iowa,
Iowa City, IA 52242, USA}
\email{thomas-branson@uiowa.edu}
\author{Andreas \v{C}ap}
\address{Fakult\"at f\"ur Mathematik, Universit\"at Wien,
Nordbergstra{\ss}e~15, A-1090 Wien,\newline\indent Austria}
\email{andreas.cap@esi.ac.at}
\author{Michael Eastwood}
\address{Department of Mathematics, University of Adelaide,
South Australia 5005}
\email{meastwoo@maths.adelaide.edu.au}
\author{A.~Rod Gover}
\address{Department of Mathematics, University of Auckland,
Private Bag 92019,\newline\indent Auckland, New Zealand}
\email{gover@math.auckland.ac.nz}
\thanks{The authors would like to thank the American Institute of Mathematics,
the Erwin Schr\"odinger Institute, the Institute for Mathematical Sciences
at the National University of Singapore, and the Banff International Research
Station for hospitality during the preparation of this article. This research
was also supported by the US NSF (Grant INT-9724781), the Austrian FWF
(Project P15747), the Australian Research Council, the Royal Society of
New Zealand (Marsden Grant 02-UOA-108), and the New Zealand Institute of
Mathematics and its Applications. The authors express their thanks to the 
referee for suggesting clarifications in the text.}
\copyrightinfo{2004}{American Mathematical Society}
\begin{abstract}
We show that a wide class of geometrically defined overdetermined semilinear
partial differential equations may be explicitly prolonged to obtain closed
systems. As a consequence, in the case of linear equations we extract sharp
bounds on the dimension of the solution space.
\end{abstract}
\renewcommand{\subjclassname}{\textup{2000} Mathematics Subject Classification}
\maketitle

\section{Introduction}
For ordinary differential equations, it is clear that the $n^{\mathrm{th}}$
order equation
$$\frac{d^n\sigma}{dx^n}=
f\Big(x,\sigma,\frac{d\sigma}{dx},\dots,\frac{d^{n-1}\sigma}{dx^{n-1}}\Big)$$
is equivalent to the system of first order equations
$$\frac{d\sigma}{dx}=\sigma_1,\,\ldots,\frac{d\sigma_k}{dx}=\sigma_{k+1},\,
\ldots,\frac{d\sigma_{n-1}}{dx}=f(x,\sigma,\sigma_1,\dots,\sigma_{n-1}).$$
This man{\oe}uvre is well-known, for example, in reducing the existence and
uniqueness of solutions to ordinary differential equations to the case of first
order equations.

For partial differential equations, however, this na\"{\i}ve man{\oe}uvre
fails. Even for overdetermined equations, it is necessary to introduce new
dependent variables for certain higher derivatives in order to achieve a first
order `closed system'---one in which all the first partial derivatives of all
the dependent variables are determined in terms of the variables themselves.
Example~\ref{ex2} below is typical in this regard---the original equation is
first order but the closed system (\ref{system12}) implicitly but necessarily
involves second derivatives of the original dependent variables~$\sigma_a$. The
introduction of new variables for unknown higher derivatives with the aim of
expressing all their derivatives as differential consequences of the original
equation is the well-known procedure of `prolongation'.

Classically, the prolongations of a semilinear differential operator
$D:E\to F$ between smooth vector bundles $E$ and $F$ on a smooth manifold $M$
are constructed from its leading symbol
$\sigma(D):\bigodot^k\Lambda^1\otimes E\to F$ where $\bigodot^k\Lambda^1$
denotes the bundle of symmetric covariant tensors on $M$ of valence~$k$. At
any point of~$M$,  denoting by $K$ the kernel of~$\sigma(D)$, one considers the
vector spaces
\begin{equation}\label{classicalV}\textstyle
V_i=(\bigodot^i\Lambda^1\otimes E)\cap(\bigodot^{i-k}\Lambda^1\otimes K)
\quad\mbox{for }i\geq k,\end{equation}
declaring the system to be of {\em finite type} if $V_i=0$ for $i$ sufficiently
large~\cite{spencer}. The solutions of a system of finite type are determined
by finitely many jets at a point. Although there is a general criterion that
$D$ be of finite type (namely, that its characteristic variety be empty
\cite[Proposition~1.7.5]{spencer}) the computation of $V_i$ presents a major
obstacle to further progress.

There are two main points to this article. Firstly, for a wide class of
geometric overdetermined partial differential equations, we explicitly
compute~$V_i$ (Lemma~\ref{injectusw} part~(4)). The direct sum
$V=\bigoplus_{i=0}^NV_i$ is a vector bundle induced by an irreducible
representation of a reductive Lie algebra so $N$ and its rank can be
immediately read off. This gives sharp bounds on the jet needed to pin down a
solution and, in the linear case, the dimension of the space of solutions. The
second point to this article is motivated by geometric considerations. We can
deal with all symbols of overdetermined invariant operators for an important
class of structures including conformal and quaternionic geometries. Motivated
by the machinery of Bernstein-Gelfand-Gelfand sequences~\cite{cd,css}, we find
a uniform procedure to perform the further steps necessary explicitly to
rewrite the equation in closed form. For the whole development, representation
theory, especially Kostant's algebraic Hodge theory~\cite{kostant} in Lie
algebra cohomology, provides the key to our method.

For readers unfamiliar with overdetermined systems, we
begin by discussing some examples. The reader should be aware, however, that
these examples are far too simple satisfactorily to illustrate the general
procedure. In fact, this is inevitable---though our algorithm is explicit, the
details in any particular case will generally be fearsome. However, for many
purposes, the details are unnecessary. For example, we may deduce without
hesitation that, on a Riemannian manifold of dimension $n\geq 3$, the space of
solutions of the partial differential equation
\begin{equation}\label{pde}\mbox{the trace-free part of }
\nabla_{(a}\nabla_b\sigma_{c)}=0\end{equation}
is finite-dimensional of dimension at most $n(n+2)(n+4)/3$. This bound is sharp
and any solution is determined by its $4$-jet at one point.

In (\ref{pde}) and throughout, we adopt Penrose's abstract index
notation~\cite{OT}. Thus, indices act as markers to specify the type of a
tensor (so $\sigma_a$ is a $1$-form whilst $\sigma^a$ would be a vector field)
and to record symmetries and contractions. Round brackets, as in~(\ref{pde}),
mean that the indices they enclose are symmetrised, square brackets
$\phi_{[ab]c}$ take the skew part, and a repeated index $\phi^a{}_{ab}$ denotes
contraction. On a Riemannian manifold, indices may be raised or lowered with
the metric in the usual way. Connections will be denoted $\nabla_a$ and on a
Riemannian manifold will usually mean the Levi-Civita connection. If $\nabla_a$
is a torsion-free connection on the tangent bundle, then its curvature tensor
$R_{ab}{}^c{}_d$ is defined by
$$(\nabla_a\nabla_b-\nabla_b\nabla_a)V^c=R_{ab}{}^c{}_dV^d.$$
In particular,
$\nabla_b\nabla_aV^b=\nabla_a\nabla_bV^b+R_{ab}V^b$,
where $R_{ab}=R_{ca}{}^c{}_b$ is the Ricci tensor.

\subsection{Two affine examples} Here we work on a smooth manifold with
torsion-free connection~$\nabla_a$.
\subsubsection{Example}\label{ex1} Consider the partial differential
equation on the function~$\sigma$:--
\begin{equation}\label{firstexample}\nabla_a\nabla_b\sigma=0.\end{equation}
If we introduce $\mu_a=\nabla_a\sigma$, then we can rewrite it as a system:--
$$\begin{array}{rcl}\nabla_a\sigma&=&\mu_a\\ \nabla_a\mu_b&=&0.\end{array}$$
\subsubsection{Example}\label{ex2} Consider the partial differential
equation on the $1$-form~$\sigma_a$:--
\begin{equation}\label{secondexample}\nabla_{(a}\sigma_{b)}=0.\end{equation}
We can rewrite it as
$$\nabla_a\sigma_b=\mu_{ab}\quad\mbox{where $\mu_{ab}$ is skew.}$$
Na\"{\i}vely differentiating this equation leads nowhere but notice that, as
differential forms, $\mu=d\sigma$ whence~$d\mu=0$. In index notation
$\nabla_{[a}\mu_{bc]}=0$ so
$$\nabla_a\mu_{bc}=\nabla_c\mu_{ba}-\nabla_b\mu_{ca}
 =\nabla_c\nabla_b\sigma_a-\nabla_b\nabla_c\sigma_a
 =R_{bc}{}^d{}_a\sigma_d.$$
Therefore, the differential equation (\ref{secondexample}) is equivalent to the
system
\begin{equation}\label{system12}
\begin{array}{rcl}\nabla_a\sigma_b&=&\mu_{ab}
\qquad\mbox{where $\mu_{ab}$ is skew}\\
\nabla_a\mu_{bc}&=&R_{bc}{}^d{}_a\sigma_d.\end{array}\end{equation}

\subsection{Two Riemannian examples}
Here we work on $n$-dimensional Riemannian manifold with metric $g_{ab}$ and
Levi-Civita connection~$\nabla_a$. We shall suppose that~$n\geq 3$.
\subsubsection{Example}\label{ex3} Consider the partial differential equation
\begin{equation}\label{thirdexample}\mbox{the trace-free part of }
\nabla_a\nabla_b\sigma=0.\end{equation}
If we introduce $\mu_a=\nabla_a\sigma$, then we can rewrite it as
$$\nabla_a\mu_b=\rho g_{ab}\quad\mbox{for some smooth function~$\rho$.}$$
Then
$$\nabla_a\rho=\nabla^b\nabla_a\mu_b=\nabla_a\nabla^b\mu_b+R_a{}^c\mu_c
 =n\nabla_a\rho+R_a{}^b\mu_b.$$
Therefore, the differential equation (\ref{thirdexample}) is equivalent to the
system
\begin{equation}\label{thirdsystem}
\begin{array}{rcl}\nabla_a\sigma&=&\mu_a\\ \nabla_a\mu_b&=&\rho g_{ab}\\
\nabla_a\rho&=&-\ffrac{1}{n-1}R_a{}^b\mu_b.\end{array}\end{equation}
\subsubsection{Example}\label{ex4}  Consider the partial differential equation
\begin{equation}\label{fourthexample}\mbox{the trace-free part of }
\nabla_{(a}\sigma_{b)}=0.\end{equation}
Even in this simple case, prolongation is already fairly involved. The details
can be omitted on first reading and the main features are described in
\S\ref{discuss} below. We can rewrite (\ref{fourthexample}) as
\begin{equation}\label{one}
\nabla_a\sigma_b=\mu_{ab}+\nu g_{ab}\quad\mbox{where $\mu_{ab}$ is skew.}
\end{equation}
Then $\nabla_{[a}\mu_{bc]}=0$, so
\begin{equation}\label{nablamu}
\begin{array}{rcl}\nabla_a\mu_{bc}=\nabla_c\mu_{ba}-\nabla_b\mu_{ca}&=&
\nabla_c(\nabla_b\sigma_a-\nu g_{ba})-\nabla_b(\nabla_c\sigma_a-\nu g_{ca})\\
&=&R_{bc}{}^d{}_a\sigma_d-g_{ab}\nabla_c\nu+g_{ac}\nabla_b\nu.
\end{array}\end{equation}
Tracing over $a$ and $b$ gives
$$\nabla^b\mu_{bc}=-R_c{}^d\sigma_d-(n-1)\nabla_c\nu.$$
Let us introduce $\rho_c=\ffrac{1}{n-1}\nabla^b\mu_{bc}$ and rearrange this
last equation as
\begin{equation}\label{two}
\textstyle\nabla_a\nu=-\rho_a-\ffrac{1}{n-1}R_a{}^b\sigma_b.\end{equation}
It may be used to eliminate $\nabla_c\nu$ from (\ref{nablamu}) to obtain
\begin{equation}\label{three}\nabla_a\mu_{bc}
=g_{ab}\rho_c-g_{ac}\rho_b+K_{abc},\end{equation}
where
\begin{equation}\label{defofK}\textstyle K_{abc}=
R_{bc}{}^d{}_a\sigma_d
+\ffrac{1}{n-1}g_{ab}R_c{}^d\sigma_d-\ffrac{1}{n-1}g_{ac}R_b{}^d\sigma_d.
\end{equation}
Notice that $K_{abc}$ is totally trace-free. Now apply $\nabla_d$ to
(\ref{three}) and skew over $a$ and $d$ to obtain
$$R_{da}{}^e{}_b\mu_{ce}-R_{da}{}^e{}_c\mu_{be}=\nabla_dK_{abc}-\nabla_aK_{dbc}
+g_{ab}\nabla_d\rho_c-g_{db}\nabla_a\rho_c
-g_{ac}\nabla_d\rho_b+g_{dc}\nabla_a\rho_b.$$
Tracing over $a$ and $b$ gives
$$R_d{}^e\mu_{ce}-R_d{}^{be}{}_c\mu_{be}=-\nabla^bK_{dbc}
+(n-2)\nabla_d\rho_c+g_{dc}\nabla^b\rho_b$$
but tracing again, over $c$ and $d$, gives $0=2(n-1)\nabla^b\rho_b$. Therefore,
\begin{equation}\label{four}\textstyle\nabla_a\rho_b=\ffrac{1}{n-2}\left(
R_a{}^c\mu_{bc}-R_a{}^{cd}{}_b\mu_{cd}-\nabla^cK_{abc}\right).\end{equation}
At this point it is clear that the system has closed: it comprises (\ref{one}),
(\ref{two}), (\ref{three}), and in (\ref{four}) one has to expand
$\nabla^cK_{abc}$ using (\ref{defofK}) and (\ref{one}).

\subsection{Discussion}\label{discuss}
In each of the examples above, we start with a linear differential operator
$D:E\to F$ between vector bundles and the conclusion is that various auxiliary
fields may be introduced so that the equation $D\sigma=0$ is equivalent to a
`closed system' in which all the first partial derivatives of all fields are
determined as linear expressions in the fields themselves. It is convenient to
regard this system as a vector bundle $\VB$ with connection~$\widetilde\nabla$.
Thus, the conclusion of Example~\ref{ex2} is that
$$\nabla_{(a}\sigma_{b)}=0\quad\mbox{if and only if}\quad
\widetilde\nabla\Sigma=0$$
where
$$\Sigma=
\left\lgroup\!\begin{array}c\sigma_b\\ \mu_{bc}\end{array}\!\right\rgroup
\quad\mbox{is a section of the vector bundle}\quad\VB=
\begin{array}c\Lambda^1\\[-2pt] \oplus\\ \Lambda^2\end{array}$$
and $\widetilde\nabla:\VB\to\Lambda^1\otimes\VB$ is the connection:--
$$\widetilde\nabla_a
\left\lgroup\!\begin{array}c\sigma_b\\ \mu_{bc}\end{array}\!\right\rgroup=
\left\lgroup\!\begin{array}c\nabla_a\sigma_b-\mu_{ab}\\
\nabla_a\mu_{bc}-R_{bc}{}^d{}_a\sigma_d\end{array}\!\right\rgroup.$$

Our examples, constructing $\Sigma$ and $\widetilde\nabla$ from $\sigma$
and~$D$, follow the well-known method of `prolongation'. Our aim in this
article, however, is to predict the form of a valid prolongation for a natural
and extensive class of examples without having to carry out the prolongation in
detail.

The conclusion of Example~\ref{ex4} is that (\ref{fourthexample}) is equivalent
to $\widetilde\nabla\Sigma=0$ where
$$\Sigma=
\left\lgroup\!\begin{array}c\sigma_b\\ \mu_{bc}\qquad\nu\\
\rho_b\end{array}\!\right\rgroup
\quad\mbox{is a section of the bundle}\quad\VB=
\begin{array}c{}\:\Lambda^1\\
\Lambda^2\bigoplus\Lambda^0\\[2pt] {}\:\Lambda^1\end{array}$$
and $\widetilde\nabla:\VB\to\Lambda^1\otimes\VB$ is an explicit connection of
the form
\begin{equation}\label{form}\widetilde\nabla
\left\lgroup\!\begin{array}c\sigma\\ \mu\qquad\nu\\
\rho\end{array}\!\right\rgroup
=\left\lgroup\!\begin{array}c\nabla\sigma-\mu-\nu\\
\nabla\mu-\rho-R\bowtie\sigma\qquad\quad\nabla\nu-\rho-R\bowtie\sigma\\
\nabla\rho-R\bowtie\mu-R\bowtie\nu-(\nabla R)\bowtie\sigma
\end{array}\!\right\rgroup,
\end{equation}
where each $\bowtie$ indicates an appropriate linear combination of
contractions of its \mbox{ingredients}.

Note that $\Sigma$ is obtained from $\sigma_a$ by application of a linear
second order differential operator, explicitly
$$\sigma_a\longmapsto
\left\lgroup\!\begin{array}c\sigma_a\\
\nabla_{[a}\sigma_{b]}\qquad\ffrac{1}{n}\nabla^a\sigma_a\\
\ffrac{1}{2(n-1)}
\left(\nabla^b\nabla_b\sigma_a-\nabla^b\nabla_a\sigma_b\right)
\end{array}\!\right\rgroup.$$

The equation (\ref{fourthexample}) is well-known. It says that the vector field
$\sigma^a$ is a conformal Killing field---its flow preserves the metric up to
scale. {From} this geometric interpretation it follows easily that the space of
solutions is bounded by $\dim{\mathfrak{so}}(n+1,1)$ since
${\mathfrak{so}}(n+1,1)$ is the conformal algebra in the flat case. This bound
is confirmed by the technique of prolongation:--
$$\rank\VB=2\rank\Lambda^1+\rank\Lambda^2+\rank\Lambda^0=
2n+\frac{n(n-1)}2+1=\frac{(n+1)(n+2)}2.$$
In~\cite{semmelmann}, Semmelmann uses this technique to establish similar
bounds on the dimension of spaces of conformal Killing forms. Specifically, he
finds an explicit connection (also having the form~(\ref{form})) on the bundle
$$\VB=\begin{array}c{}\:\Lambda^p\\
\Lambda^{p+1}\bigoplus\Lambda^{p-1}\\[2pt] {}\:\Lambda^p\end{array}
\quad\mbox{with rank }
\left(\!\begin{array}cn+2\\ p+1\end{array}\!\right)$$
so that conformal Killing $p$-forms are equivalent to parallel sections of this
bundle. The general procedure, to be explained in this article, includes this
case and many more besides.

The corresponding bound for Example~\ref{ex2} is
$$\rank\Lambda^1+\rank\Lambda^2=n+\frac{n(n-1)}2=\frac{n(n+1)}2.$$
It was pointed out to us by Dan Fox that this is precisely the bound
investigated by Eisenhart in~\cite{eisenhart}.

\subsection{Semilinear variants}
Each of the examples discussed so far persists in a semilinear form. Thus,
Example~\ref{ex1} may be modified as
$$\nabla_a\nabla_b\sigma=f_{ab}(x,\sigma,\nabla_c\sigma)$$
where $f_{ab}$ depends smoothly on its arguments and takes values in symmetric
2-tensors. Evidently, this equation is equivalent to the system
$$\begin{array}{rcl}\nabla_a\sigma&=&\mu_a\\
\nabla_a\mu_b&=&f_{ab}(x,\sigma,\mu_a).\end{array}$$

Example~\ref{ex2} may be modified as
\begin{equation}\label{quasi2}\nabla_{(a}\sigma_{b)}=f_{ab}(x,\sigma_c).
\end{equation}
The only difficulty in following previous reasoning is that one must be careful
as to the meaning of $\nabla_cf_{ab}(x,\sigma_d)$. As it arises, $\sigma_d$ is
a function of $x$ and so $f_{ab}$ may be regarded as a tensor on the manifold
and $\nabla_cf_{ab}$ as the usual covariant derivative. On the other hand, we
may fix~$\sigma_d$, regard $f_{ab}(x,\sigma_d)$ as a function of its first
argument, and then take its covariant derivative. We shall use the notation
$\partial_cf_{ab}$ for the result of this point of view. There is also the
partial derivative obtained by fixing $x$ and differentiating with respect
to~$\sigma_d$: let us write $\delta^d=\partial/\partial\sigma_d$. Then, by the
chain rule,
$$\nabla_cf_{ab}=\partial_cf_{ab}+(\delta^df_{ab})\nabla_c\sigma_d,$$
often referred to as expressing `total derivative' in terms of `partial
derivative'. The result of following previous reasoning is that (\ref{quasi2})
is equivalent to the system
$$\begin{array}{rcl}\nabla_a\sigma_b&=&\mu_{ab}+f_{ab}\\
\nabla_a\mu_{bc}&=&R_{bc}{}^d{}_a\sigma_d
+2\partial_{[b}f_{c]a}-2(\delta^df_{a[b})\mu_{c]d}-2(\delta^df_{a[b})f_{c]d},
\end{array}$$
where $\mu_{ab}$ is skew. As a typical nonlinear variant therefore,
$$\nabla_{(a}\sigma_{b)}=\sigma_a\sigma_b+S_{ab}$$
for an arbitrary given symmetric tensor $S_{ab}$ is equivalent to the closed
system
$$\begin{array}{rcl}\nabla_a\sigma_b&=&\mu_{ab}+\sigma_a\sigma_b+S_{ab}\\
\nabla_a\mu_{bc}&=&R_{bc}{}^d{}_a\sigma_d+2\nabla_{[b}S_{c]a}
-2\sigma_{[b}\mu_{c]a}+2\sigma_a\mu_{bc}+2S_{a[b}\sigma_{c]}.
\end{array}$$

The general semilinear variant of Example~\ref{ex3} is
$$\textstyle\nabla_a\nabla_b\sigma-\ffrac1ng_{ab}\nabla^c\nabla_c\sigma=
f_{ab}(x,\sigma,\nabla_c\sigma),$$
where $f_{ab}(x,\sigma,\sigma_c)$ is symmetric and trace-free. If we write
$\delta=\partial/\partial\sigma$, then the chain rule for total derivative in
terms of partial derivative is
$$\nabla_cf_{ab}=\partial_cf_{ab}+(\delta f_{ab})\nabla_c\sigma
+(\delta^df_{ab})\nabla_c\sigma_d$$
and the closed system generalising (\ref{thirdsystem}) is
$$\begin{array}{rcl}\nabla_a\sigma&=&\mu_a\\
\nabla_a\mu_b&=&\rho g_{ab}+f_{ab}\\
\nabla_a\rho&=&-\ffrac{1}{n-1}R_a{}^b\mu_b
+\ffrac{1}{n-1}(\partial^bf_{ab}+(\delta f_{ab})\mu^b
+(\delta^bf_{ab})\rho+(\delta^df_{ab})f^b{}_d).\end{array}$$

A particular semilinear variant of Example~\ref{ex4} is
$$\textstyle\mbox{the trace-free part of }
(\nabla_{(a}\sigma_{b)}+\sigma_a\sigma_b+\ffrac{1}{n-2}R_{ab})=0,$$
where $R_{ab}$ is the Ricci tensor. It is the Einstein-Weyl equation and the
corresponding closed system is derived in~\cite{e-w} by {\em ad hoc\/} methods.

\section{Formulation of the main results}\label{formulation}
Firstly, some generalities on differential operators. As detailed
in~\cite{spencer}, to every smooth vector bundle $E$ on a smooth manifold $M$
there are the canonically associated jet bundles $J^kE$ on $M$ and short exact
sequences of homomorphisms of vector bundles
$$\textstyle0\to\bigodot^k\Lambda^1\otimes E\to J^kE\to J^{k-1}E\to 0,$$
where $\bigodot^k\Lambda^1$ denotes the $k^{\mathrm{th}}$ symmetric tensor
power of~$\Lambda^1$. A $k^{\mathrm{th}}$ order linear differential operator
$D:E\to F$ between vector bundles $E$ and $F$ is equivalent to a homomorphism
of vector bundles $J^kE\to F$ and the symbol $\sigma(D)$ of $D$ is defined as
the composition
\begin{equation}\label{symbol}
\textstyle\bigodot^k\Lambda^1\otimes E\hookrightarrow J^kE\to F.\end{equation}
A differential operator of the form $D_1+D_2$ where $D_1$ is $k^{\mathrm{th}}$
order linear and $D_2$ is $(k-1)^{\mathrm{st}}$ order is called semilinear and
its symbol is defined to be~$\sigma(D_1)$.

If we now return to the semilinear variants of our affine examples, we see that
the form of the equation is independent of the connection.
Equation~(\ref{quasi2}), for example, says that we have a first order
semilinear operator $\Lambda^1\to\bigodot^2\Lambda^1$ whose symbol
$$\textstyle\Lambda^1\otimes\Lambda^1\longrightarrow\bigodot^2\Lambda^1$$
is taking the symmetric part. In particular, a change of torsion-free
connection in~(\ref{secondexample}) is covered by
$$\nabla_{(a}\sigma_{b)}=\Gamma_{ab}{}^c\sigma_c$$
as a special case of~(\ref{quasi2}).

To formulate the semilinear equations on a smooth manifold $M$ to which our
prolongation procedure will apply, let us regard the tangent bundle as
tautologically associated to the frame bundle under the standard representation
of ${\mathrm{GL}}(n,{\mathbb{R}})$ on~${\mathbb{R}}^n$. Then, an irreducible
tensor bundle on $M$ is, by definition, a bundle associated to the frame bundle
under an irreducible representation of ${\mathrm{GL}}(n,{\mathbb{R}})$. By
basic representation theory, any tensor bundle decomposes into a direct sum of
irreducible tensor bundles. In fact, for technical reasons, let us fix a volume
form on~$M$. This reduces the structure group of the frame bundle to
${\mathrm{SL}}(n,{\mathbb{R}})$ and allows us to use the usual theory of
weights to specify an irreducible tensor bundle. Following~\cite{beastwood},
the irreducible representations are in one-to-one correspondence with
attachments of non-negative integers to the nodes of the Dynkin diagram of
${\mathfrak{sl}}(n,{\mathbb{R}})$. 
These numbers represent the coefficients in the expansion of the highest 
weight of the dual representation (or equivalently the negative of the lowest 
weight of the given representation) as a linear combination of fundamental 
weights. Each coefficient is placed over the node representing the simple 
root that is dual to the fundamental weight. 
Combining these viewpoints, the tangent bundle is
$$\oooospecial{0}{0}{0}{1}$$
and so on:--
$$\textstyle\Lambda^1=\oooo{1}{0}{0}{0}\qquad
\Lambda^2=\oooo{0}{1}{0}{0}\qquad
\bigodot^k\Lambda^1=\oooo{k}{0}{0}{0}.$$
For any irreducible tensor bundle~$E$, the tensor product
$\bigodot^k\Lambda^1\otimes E$ decomposes into irreducibles amongst which, the
one with highest weight obtained by adding the highest weights of
$\bigodot^k\Lambda^1$ and~$E$, occurs with multiplicity one. This is the
Cartan product \cite{dynkin} and we shall denote it
$\bigodot^k\Lambda^1\cart E$. In the notation just established,
$$\textstyle E=\oooo{a}{b}{c}{d}\quad\Rightarrow\quad
\bigodot^k\Lambda^1\cart E=\oooo{k+a}{b}{c}{d}.$$
There is a canonical projection
$\bigodot^k\Lambda^1\otimes E\to\bigodot^k\Lambda^1\cart E$, which we shall
also refer to as the Cartan product.

Now we are in a position to state a special case of our main theorem:--
\begin{theorem}\label{specialaffineversion}
Suppose $M$ is a smooth manifold equipped with a volume form. Let $E$ be an
irreducible tensor bundle on $M$ and $F=\bigodot^k\Lambda^1\cart E$. Suppose
$D:E\to F$ is a $k^{\mathrm{th}}$-order semilinear differential operator whose
symbol
$$\textstyle
\sigma(D):\bigodot^k\Lambda^1\otimes E\to F=\bigodot^k\Lambda^1\cart E$$
is given by the Cartan product. Then, there is a vector bundle $\VB$ and, for
every choice of volume-preserving connection $\nabla$ on the tangent bundle, a
canonically associated connection $\widetilde\nabla:\VB\to\Lambda^1\otimes\VB$
on $\VB$ so that there is a bijection
\begin{equation}\label{isomorphism}
\{\sigma\in\Gamma(E)\mathrm{\ s.t.\ }D\sigma=0\}\cong
\{\Sigma\in\Gamma(\VB)\mathrm{\ s.t.\ }
\widetilde\nabla\Sigma+\Phi(\Sigma)=0\},\end{equation}
where $\Phi:\VB\to\Lambda^1\otimes\VB$ is a fibre-preserving
mapping canonically constructed from~$D$. If $D$ is linear, then so is~$\Phi$.
{From} left to right {\rm (\ref{isomorphism})} is implemented by an
$N^{\mathrm{th}}$ order linear differential operator where $N$ is easily
computable from $E$ and~$k$.
\end{theorem}
\noindent We should emphasise that the bundle $\VB$ is completely determined by
$E$ and~$k$. The connection $\widetilde\nabla$ on $\VB$ is then determined by a
choice of affine connection on~$M$. Finally, the fibre-preserving mapping
$\Phi$ is determined by~$D$.

In particular, $\VB$ is defined as follows. Let us embed
${\mathrm{SL}}(n,{\mathbb{R}})\hookrightarrow{\mathrm{SL}}(n+1,{\mathbb{R}})$
by $${\mathrm{SL}}(n,{\mathbb{R}})\ni A\longmapsto
\mbox{\small$\left\lgroup\!\!\begin{array}{cc}
1&0\\ 0&A\end{array}\!\!\right\rgroup$}
\in{\mathrm{SL}}(n+1,{\mathbb{R}}).$$
Corresponding to this embedding, the Dynkin diagram of 
${\mathfrak{sl}}(n+1,{\mathbb{R}})$ is obtained from the Dynkin diagram of 
${\mathfrak{sl}}(n,{\mathbb{R}})$ by adding a node on the left. Let us denote 
the fundamental weight of ${\mathfrak{sl}}(n+1,{\mathbb{R}})$ corresponding 
to the additional simple root by $\omega_0$. Any representation of 
${\mathrm{SL}}(n+1,{\mathbb{R}})$ restricts to a
representation of ${\mathrm{SL}}(n,{\mathbb{R}})$ and hence gives rise to an
associated vector bundle on~$M$. Using these two facts, given  
$E=\oooo{a}{b}{c}{d}$ and $k\geq 1$, we define 
$\VB:=\ooooo{k-1}{a}{b}{c}{d}$, and it turns out that 
$N=k-1+a+b+c+\cdots+d$. More explicitly, if $E$ is associated to the dual of 
the irreducible representation of ${\mathrm{SL}}(n,{\mathbb{R}})$ with 
highest weight $\lambda$, then we consider the irreducible representation of 
${\mathrm{SL}}(n+1,{\mathbb{R}})$ with 
highest weight $(k-1)\omega_0+\lambda$, restrict its dual to 
${\mathrm{SL}}(n,{\mathbb{R}})$ and let $V$ be the associated vector bundle. 
When restricted to ${\mathrm{SL}}(n,{\mathbb{R}})$, an irreducible 
representation of ${\mathrm{SL}}(n+1,{\mathbb{R}})$ splits into a direct sum 
of irreducible representations of ${\mathrm{SL}}(n,{\mathbb{R}})$. 
Correspondingly, we obtain a splitting
\begin{equation}\label{affinesplit}
\VB=\ooooo{k-1}{a}{b}{c}{d}=\oooo{a}{b}{c}{d}\oplus\quad\cdots.\end{equation}
The representation corresponding to the first summand, which is (isomorphic to)
$E$, can be described as the ${\mathrm{SL}}(n,{\mathbb{R}})$--invariant 
subspace generated by a vector of lowest weight. In particular, there is a
canonically defined surjection $\pi:\VB\to E$ and it is $\sigma=\pi\circ\Sigma$
that induces the isomorphism~(\ref{isomorphism}) from right to left.

In the situation of Example \ref{ex1}, $E$ corresponds to the trivial 
representation and $k=2$. Thus we obtain $\VB=\ooooo{1}{0}{0}{0}{0}$. This 
corresponds to the representation ${\mathbb{R}}^{(n+1)*}$, which restricted 
to ${\mathrm{SL}}(n,{\mathbb{R}})$ splits as ${\mathbb{R}}\oplus
{\mathbb{R}}^{n*}$. Hence we obtain $\VB={\mathbb{R}}\oplus\Lambda^1$ and 
$N=1$.

For Example \ref{ex2}, we have $E=\Lambda^1$ and $k=1$, which implies 
$\VB=\ooooo{0}{1}{0}{0}{0}$. The corresponding representation 
$\Lambda^2{\mathbb{R}}^{(n+1)*}$ splits as 
${\mathbb{R}}^{n*}\oplus\Lambda^2{\mathbb{R}}^{n*}$, 
so $\VB=\Lambda^1\oplus\Lambda^2$ and again $N=1$. 

For the Riemannian version of Theorem~\ref{specialaffineversion} we simply
replace the embedding of Lie groups
${\mathrm{SL}}(n,{\mathbb{R}})\hookrightarrow{\mathrm{SL}}(n+1,{\mathbb{R}})$
by the embedding ${\mathrm{SO}}(n)\hookrightarrow{\mathrm{SO}}(n+1,1)$:--
$${\mathrm{SO}}(n)\ni A\longmapsto
\mbox{\small$\left\lgroup\!\!\begin{array}{ccc}
1&0&0\\ 0&A&0\\ 0&0&1\end{array}\!\!\right\rgroup$}
\in{\mathrm{SO}}(n+1,1),$$
where ${\mathrm{SO}}(n+1,1)$ is realised as preserving the quadratic form
$2x_0x_{n+1}+\sum_{i=1}^nx_i{}^2$.
There is a corresponding inclusion of Dynkin diagrams:--
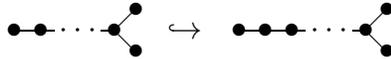
$$\begin{picture}(60,16)(10,0)
\put(12,8){\makebox(0,0){$\bullet$}}
\put(22,8){\makebox(0,0){$\bullet$}}
\put(12,8){\line(1,0){15}}
\put(50,8){\line(-1,0){5}}
\put(49,7){\line(1,1){10}}
\put(49,9){\line(1,-1){10}}
\put(37,8){\makebox(0,0){$\cdots$}}
\put(50,8){\makebox(0,0){$\bullet$}}
\put(58,16){\makebox(0,0){$\bullet$}}
\put(58,0){\makebox(0,0){$\bullet$}}
\end{picture}
\raisebox{5pt}{$\hookrightarrow$}\quad
\begin{picture}(70,16)
\put(2,8){\makebox(0,0){$\bullet$}}
\put(12,8){\makebox(0,0){$\bullet$}}
\put(22,8){\makebox(0,0){$\bullet$}}
\put(2,8){\line(1,0){25}}
\put(50,8){\line(-1,0){5}}
\put(49,7){\line(1,1){10}}
\put(49,9){\line(1,-1){10}}
\put(37,8){\makebox(0,0){$\cdots$}}
\put(50,8){\makebox(0,0){$\bullet$}}
\put(58,16){\makebox(0,0){$\bullet$}}
\put(58,0){\makebox(0,0){$\bullet$}}
\end{picture}$$
if $n$ is even and
$$\begin{picture}(60,16)(10,0)
\put(12,8){\makebox(0,0){$\bullet$}}
\put(22,8){\makebox(0,0){$\bullet$}}
\put(12,8){\line(1,0){15}}
\put(50,8){\line(-1,0){5}}
\put(37,8){\makebox(0,0){$\cdots$}}
\put(55,8){\makebox(0,0){$\rangle$}}
\put(50,8){\makebox(0,0){$\bullet$}}
\put(60,8){\makebox(0,0){$\bullet$}}
\put(50,9){\line(1,0){10}}
\put(50,7){\line(1,0){10}}
\end{picture}
\raisebox{5pt}{$\hookrightarrow$}\quad
\begin{picture}(70,16)
\put(2,8){\makebox(0,0){$\bullet$}}
\put(12,8){\makebox(0,0){$\bullet$}}
\put(22,8){\makebox(0,0){$\bullet$}}
\put(2,8){\line(1,0){25}}
\put(50,8){\line(-1,0){5}}
\put(37,8){\makebox(0,0){$\cdots$}}
\put(55,8){\makebox(0,0){$\rangle$}}
\put(50,8){\makebox(0,0){$\bullet$}}
\put(60,8){\makebox(0,0){$\bullet$}}
\put(50,9){\line(1,0){10}}
\put(50,7){\line(1,0){10}}
\end{picture}$$
if $n$ is odd. The irreducible tensor bundles on an oriented Riemannian
manifold are associated to irreducible representations of ${\mathrm{SO}}(n)$.
On an oriented spin manifold, we should use
${\mathrm{Spin}}(n)\hookrightarrow{\mathrm{Spin}}(n+1,1)$ instead and
there are irreducible spinor bundles too, associated to irreducible spin
representations. The Riemannian version of Theorem~\ref{specialaffineversion}
is obtained by taking $F=\bigodot_\circ^k\Lambda^1\cart E$ where
$\bigodot_\circ$ denotes trace-free symmetric product. Thus, if $n$ is odd for
example, then
$$E=\raisebox{-6pt}{\begin{picture}(60,20)(10,0)
\put(12,13){\makebox(0,0)[b]{$\scriptstyle a$}}
\put(22,13){\makebox(0,0)[b]{$\scriptstyle b$}}
\put(50,13){\makebox(0,0)[b]{$\scriptstyle c$}}
\put(60,13){\makebox(0,0)[b]{$\scriptstyle d$}}
\put(12,8){\makebox(0,0){$\bullet$}}
\put(22,8){\makebox(0,0){$\bullet$}}
\put(12,8){\line(1,0){15}}
\put(50,8){\line(-1,0){5}}
\put(37,8){\makebox(0,0){$\cdots$}}
\put(55,8){\makebox(0,0){$\rangle$}}
\put(50,8){\makebox(0,0){$\bullet$}}
\put(60,8){\makebox(0,0){$\bullet$}}
\put(50,9){\line(1,0){10}}
\put(50,7){\line(1,0){10}}
\end{picture}}
\quad\Rightarrow\quad
\mbox{\Large$\Big\{$}\begin{array}l\raisebox{6pt}{$F={}$}
\begin{picture}(70,20)
\put(7,8){\makebox(0,0){$\bullet$}}
\put(7,13){\makebox(0,0)[b]{$\scriptstyle k+a$}}
\put(22,13){\makebox(0,0)[b]{$\scriptstyle b$}}
\put(50,13){\makebox(0,0)[b]{$\scriptstyle c$}}
\put(60,13){\makebox(0,0)[b]{$\scriptstyle d$}}
\put(22,8){\makebox(0,0){$\bullet$}}
\put(7,8){\line(1,0){20}}
\put(50,8){\line(-1,0){5}}
\put(37,8){\makebox(0,0){$\cdots$}}
\put(55,8){\makebox(0,0){$\rangle$}}
\put(50,8){\makebox(0,0){$\bullet$}}
\put(60,8){\makebox(0,0){$\bullet$}}
\put(50,9){\line(1,0){10}}
\put(50,7){\line(1,0){10}}
\end{picture}
\raisebox{6pt}{\quad$\VB={}$}
\begin{picture}(80,20)(-10,0)
\put(-3,8){\makebox(0,0){$\bullet$}}
\put(-3,13){\makebox(0,0)[b]{$\scriptstyle k-1$}}
\put(12,13){\makebox(0,0)[b]{$\scriptstyle a$}}
\put(22,13){\makebox(0,0)[b]{$\scriptstyle b$}}
\put(50,13){\makebox(0,0)[b]{$\scriptstyle c$}}
\put(60,13){\makebox(0,0)[b]{$\scriptstyle d$}}
\put(12,8){\makebox(0,0){$\bullet$}}
\put(22,8){\makebox(0,0){$\bullet$}}
\put(-3,8){\line(1,0){30}}
\put(50,8){\line(-1,0){5}}
\put(37,8){\makebox(0,0){$\cdots$}}
\put(55,8){\makebox(0,0){$\rangle$}}
\put(50,8){\makebox(0,0){$\bullet$}}
\put(60,8){\makebox(0,0){$\bullet$}}
\put(50,9){\line(1,0){10}}
\put(50,7){\line(1,0){10}}
\end{picture}\\
N=2(k-1+a+b+\cdots +c)+d\end{array}.$$
With these replacements, the Riemannian statement is almost identical. The only
significant difference is that we may as well use the Levi-Civita connection in
the construction of~$\widetilde\nabla$, which thereby becomes canonical.

\subsection{Other geometries}\label{others}
Though the affine and Riemannian cases are perhaps the most significant, there
is a more general formulation in terms of certain $G$-structures, which
provides a uniform approach and whose proof is no more difficult. It is this
approach that we shall adopt for the remainder of this article.

Let $G$ be a Lie group whose Lie algebra ${\mathfrak{g}}$ is $|1|$-graded
semisimple:--
$${\mathfrak{g}}=
{\mathfrak{g}}_{-1}\oplus{\mathfrak{g}}_0\oplus{\mathfrak{g}}_1$$
as, for example, discussed in~\cite{baston,css,ochiai}. Let $G_0\subset G$ be
the subgroup consisting of those elements whose adjoint action on
${\mathfrak{g}}$ preserves the grading. Its Lie algebra is~${\mathfrak{g}}_0$.
Let $G_0'$ be a subgroup of $G_0$ whose Lie algebra is
$[{\mathfrak{g}}_0,{\mathfrak{g}}_0]$. It is semisimple and the adjoint action
makes ${\mathfrak{g}}_{-1}$ into a $G_0'$-module. We shall suppose that $M$ is
a smooth manifold endowed with a first order $G_0'$-structure. More
specifically, $M$ should have the same dimension as ${\mathfrak{g}}_{-1}$ and
the frame bundle should be reduced under
$G_0'\to{\mathrm{GL}}({\mathfrak{g}}_{-1})$. 
If $G={\mathrm{SL}}(n+1,{\mathbb{R}})$, there is a $|1|$-grading on
${\mathfrak{g}}={\mathfrak{sl}}(n+1,{\mathbb{R}})$ so that 
$G_0'={\mathrm{SL}}(n,{\mathbb{R}})$, included into $G$ as in the discussion
after Theorem~\ref{specialaffineversion}. This leads to the standard 
representation of $G_0'$ on ${\mathfrak{g}}_{-1}\cong{\mathbb{R}}^n$, so 
the corresponding geometries are $n$-manifolds endowed with a volume form. For
$G={\mathrm{SO}}(n+1,1)$, we may arrange a $|1|$-grading so that 
$G_0'\hookrightarrow G$ becomes the inclusion of ${\mathrm{SO}}(n)$ described
above, and the corresponding geometries are oriented Riemannian $n$-manifolds.

For $M$ endowed with a $G_0'$-structure, as above, we may consider vector
bundles on $M$ induced from irreducible representations of~$G_0'$. If
${\mathbb{E}}$ is such a representation, we shall write $E$ for the
corresponding vector bundle. In particular, the adjoint action of $G_0'$ on
${\mathfrak{g}}_{-1}$ is irreducible and induces the tangent bundle. The
Killing form on ${\mathfrak{g}}$ canonically identifies ${\mathfrak{g}}_{-1}^*$
with ${\mathfrak{g}}_1$ as $G_0$-modules. Therefore, the $G_0'$-module
${\mathfrak{g}}_1$ gives rise to the cotangent bundle $\Lambda^1$ on~$M$. It is
convenient to write $\Cart^k\Lambda^1\cart E$ for the vector bundle associated
to the Cartan product $\Cart^k{\mathfrak{g}}_1\cart{\mathbb{E}}$.

A principal $G_0'$-connection gives rise to connections on all the associated
vector bundles~$E$. Conversely, because the $G_0'$-action on
${\mathfrak{g}}_{-1}$ is infinitesimally effective, a connection on the tangent
bundle compatible with the $G_0'$-structure, gives rise to a principal
connection. Here is the general statement extending
Theorem~\ref{specialaffineversion}:--
\begin{theorem}\label{generalversion}
Let $M$ be a manifold with $G_0'$-structure as above. Suppose $E$ is a vector
bundle on $M$ induced from an irreducible representation of~$G_0'$ and fix
$k\geq 1$. Then there is a vector bundle $\VB$ explicitly constructed from $E$
and~$k$ and, for every choice of $G_0'$-compatible connection $\nabla$ on the
tangent bundle, a canonically associated connection
$\widetilde\nabla:\VB\to\Lambda^1\otimes\VB$ on $\VB$ with the following
property. For every $k^{\mathrm{th}}$-order semilinear differential operator
$D:E\to F=\Cart^k\Lambda^1\cart E$ whose symbol
$$\textstyle
\sigma(D):\bigodot^k\Lambda^1\otimes E\to F=\Cart^k\Lambda^1\cart E$$
is the Cartan product, we have a bijection
\begin{equation}\label{genisomorphism}
\{\sigma\in\Gamma(E)\mathrm{\ s.t.\ }D\sigma=0\}\cong
\{\Sigma\in\Gamma(\VB)\mathrm{\ s.t.\ }
\widetilde\nabla\Sigma+\Phi(\Sigma)=0\}\end{equation}
(implemented by an $N^{\mathrm{th}}$ order linear differential operator in one
direction and
the natural projection
in the other), where $\Phi:\VB\to\Lambda^1\otimes\VB$
is a fibre-preserving mapping canonically constructed from~$D$. If $D$ is
linear, then so is~$\Phi$.
\end{theorem}
\noindent The proof will occupy~\S\ref{proof} but there is a useful and
immediate corollary:--
\begin{corollary}\label{dimbound} Any solution of $D\sigma=0$ is
determined by its $N$-jet. If $D:E\to F$ is linear, then the dimension of the
space of solutions of $D\sigma=0$ is bounded by $\rank\VB$.
\end{corollary}
\begin{proof}When $D$ is linear $\Phi$ is a homomorphism and so
$\widetilde\nabla+\Phi$ is a connection on~$\VB$. According to
(\ref{genisomorphism}), we seek parallel section of $\VB$ with respect to this
connection.\end{proof}
\noindent As in the affine and Riemannian cases, the bundle $\VB$ is induced
from an irreducible representation ${\mathbb{V}}$ of~$G$. Hence,
$\rank\VB=\dim{\mathbb{V}}$ and~$N$, which is related to the decomposition of
${\mathbb{V}}$ as a $G_0'$-module, can be computed by standard
tools from representation theory~\cite{fultonandharris,humphreys}.
For example, the ${\mathrm{SO}}(n+1,1)$-module $\son{1}{1}{0}{0}$ has
dimension $n(n+2)(n+4)/3$ and has~$N=4$, the stated bounds for (\ref{pde}).

Sharpness of both bounds is observed in the remarks at the end of this article.

\section{Algebraic ingredients}\label{algebra}
We shall need some results from Lie algebra cohomology. Specifically, what we
need is a special case of Kostant's algebraic Hodge theory~\cite{kostant}. In
this section, we state what we need. Proofs may be found in~\cite{kostant}. A
more general exposition in a similar context may be found in~\cite{cd}.

The setting is a $|1|$-graded Lie algebra ${\mathfrak{g}}$ corresponding to a
semisimple Lie group~$G$, as discussed in \S\ref{others}. Recall that $G_0'$ is
the semisimple part of $G_0$, which is itself a subgroup of $G$ with Lie
algebra~${\mathfrak{g}}_0$. Let ${\mathbb{V}}$ be an irreducible
representation of~$G$. We define a complex of $G_0$-modules
\begin{equation}\label{differential}\begin{array}{ccccccccc}
0&\longrightarrow&{\mathbb{V}}&\stackrel{\partial}{\longrightarrow}&
{\mathfrak{g}}_1\otimes{\mathbb{V}}&\stackrel{\partial}{\longrightarrow}&
\Lambda^2{\mathfrak{g}}_1\otimes{\mathbb{V}}&
\stackrel{\partial}{\longrightarrow}&\cdots\\
&&\|&&\|&&\|\\
0&\longrightarrow&{\mathbb{V}}&\stackrel{\partial}{\longrightarrow}&
\Hom({\mathfrak{g}}_{-1},{\mathbb{V}})&\stackrel{\partial}{\longrightarrow}&
\Hom(\Lambda^2{\mathfrak{g}}_{-1},{\mathbb{V}})
&\stackrel{\partial}{\longrightarrow}&\cdots
\end{array}\end{equation}
where the vertical identifications are by means of the Killing form and
$$\partial:\Hom(\Lambda^p{\mathfrak{g}}_{-1},{\mathbb{V}})\longrightarrow
\Hom(\Lambda^{p+1}{\mathfrak{g}}_{-1},{\mathbb{V}})$$
is defined by
$$\partial\phi(X_0,\dots,X_p)=
\sum_{i=0}^p(-1)^iX_i\phi(X_0,\dots,\widehat{X_i},\dots,X_p).$$
Since ${\mathfrak{g}}_{-1}$ is Abelian, it is easily verified that
$\partial^2=0$ and we define the Lie algebra cohomology
$$H^p({\mathfrak{g}}_{-1},{\mathbb{V}})=
\frac{\ker\partial:\Lambda^p{\mathfrak{g}}_1\otimes{\mathbb{V}}
\longrightarrow
\Lambda^{p+1}{\mathfrak{g}}_1\otimes{\mathbb{V}}}
{\mbox{$\im\partial:\Lambda^{p-1}{\mathfrak{g}}_1\otimes{\mathbb{V}}
\longrightarrow
\Lambda^p{\mathfrak{g}}_1\otimes{\mathbb{V}}$}}.$$
Since $\partial$ is a homomorphism of $G_0$-modules,
$H^p({\mathfrak{g}}_{-1},{\mathbb{V}})$ is a $G_0$-module. There is also a
codifferential
\begin{equation}\label{codiff}
0\longleftarrow{\mathbb{V}}\stackrel{\partial^*}{\longleftarrow}
{\mathfrak{g}}_1\otimes{\mathbb{V}}\stackrel{\partial^*}{\longleftarrow}
\Lambda^2{\mathfrak{g}}_1\otimes{\mathbb{V}}
\stackrel{\partial^*}{\longleftarrow}\cdots\end{equation}
defined by
$$\partial^*(Z_0\wedge\dots\wedge Z_p\otimes v)=
\sum_{i=0}^p(-1)^{i+1}
Z_0\wedge\dots\wedge\widehat{Z_i}\wedge\dots\wedge Z_p\otimes Z_iv.$$
It is also $G_0$-equivariant and satisfies $\partial^*{}^2=0$. There is a
`Hodge decomposition':--
\begin{equation}\label{hodge}\Lambda^p{\mathfrak{g}}_1\otimes{\mathbb{V}}=
\im(\partial)\oplus(\ker(\partial)\cap\ker(\partial^*))\oplus\im(\partial^*)
\end{equation}
and, in particular, a canonical isomorphism
$$H^p({\mathfrak{g}}_{-1},{\mathbb{V}})\cong\ker(\partial)\cap\ker(\partial^*)
\quad\mbox{on }\Lambda^p{\mathfrak{g}}_1\otimes{\mathbb{V}}.$$
The differential $\partial$ is seen more clearly in the Hodge decomposition
$$\begin{array}{rcc}\Lambda^p{\mathfrak{g}}_1\otimes{\mathbb{V}}&=&
\ker(\partial)\oplus\im(\partial^*)\\
\downarrow\phantom{\;{\mathbb{V}}}&&\mbox{\large$\swarrow\;{}$}\\
\Lambda^{p+1}{\mathfrak{g}}_1\otimes{\mathbb{V}}&=&
\im(\partial)\oplus\ker(\partial^*)\end{array}$$
as an isomorphism $\partial:\im(\partial^*)\to\im(\partial)$. Its inverse is
not necessarily~$\partial^*$. Instead, we may define $\delta^*$ to be this
inverse on $\im(\partial)$ and to annihilate $\ker(\partial^*)$. We obtain a
new $G_0$-equivariant codifferential defining the same Hodge decomposition as
does $\partial^*$ but with the congenial feature that
\begin{equation}\label{congenial}
\delta^*\partial={\rm{id}}\mbox{ on }\im(\delta^*)=\im(\partial^*)
\quad\mbox{and}\quad
\partial\delta^*={\rm{id}}\mbox{ on }\im(\partial).\end{equation}
Now let us be more specific about the representation~${\mathbb{V}}$.
The description of $|1|$--gradings is well known: for an appropriate
choice of a Cartan subalgebra for the complexification of
${\mathfrak{g}}$ there is a distinguished simple root $\alpha_0$. This
has the property that a root space lies in the complexification of
${\mathfrak{g}}_j$ ($j=-1,0,1$) if and only if $j$ is the coefficient
of $\alpha_0$ in the expansion of the given root into simple roots. In
particular, the Dynkin diagram of ${\mathfrak{g}}_0'$ is obtained by
removing in the Dynkin diagram of ${\mathfrak{g}}$ the node
representing $\alpha_0$ and all edges connected to that node. In the
affine and Riemannian cases previously discussed this was the leftmost node.  
Let $\omega_0$ denote the fundamental weight corresponding to~$\alpha_0$.
Starting with an irreducible representation ${\mathbb{E}}$ of~$G_0'$,
we may add $(k-1)\omega_0$ to the highest weight of ${\mathbb{E}}^*$, and 
define ${\mathbb{V}}$ as the dual of the irreducible representation of 
$G$ with that highest weight. 

The subalgebra ${\mathfrak{g}}_1\subset{\mathfrak{g}}$ is the nilradical of 
the parabolic ${\mathfrak{g}}_0\oplus{\mathfrak{g}}_1$, so 
Kostant's version of the Bott-Borel-Weil Theorem, see \cite{kostant}, 
describes the cohomology of ${\mathfrak{g}}_1$ with coefficients in an 
irreducible representation of ${\mathfrak{g}}$. It also follows from Kostant's
theory that $H^*({\mathfrak{g}}_1,{\mathbb{V}}^*)$ is dual (as a 
representation of ${\mathfrak{g}}_0$) to 
$H^*({\mathfrak{g}}_{-1},{\mathbb{V}})$. Since we use highest weights of dual
representations as labels, we can directly apply Kostant's algorithm. This 
describes the highest weights of irreducible components in the cohomology in 
terms of the actions of the elements of a subset $W^{\mathfrak{p}}$ of the 
Weyl group of ${\mathfrak{g}}$.
  
In particular, $H^0({\mathfrak{g}}_1,{\mathbb{V}}^*)$ is the irreducible 
representation of ${\mathfrak{g}}_0'$ whose highest weight is the restriction 
of the highest weight of ${\mathbb{V}}^*$, whence 
\begin{equation}\label{BBW0}
H^0({\mathfrak{g}}_{-1},{\mathbb{V}})={\mathbb{E}}.
\end{equation}
In particular, note that ${\mathbb{E}}$ has acquired the structure of a
$G_0$-module.

To deal with the first cohomology, we have to consider elements of the Weyl 
group which have length one, i.e.~are reflections corresponding to simple 
roots. The only simple reflection which lies in $W^{\mathfrak{p}}$ is the one 
corresponding to $\alpha_0$. This means that 
$H^1({\mathfrak{g}}_1,{\mathbb{V}}^*)$ is an irreducible representation of 
${\mathfrak{g}}_0'$, and its highest weight is obtained from the highest 
weight $\lambda$ of ${\mathbb{V}}^*$ by subtracting $(\ell+1)\alpha_0$, where 
$\ell$ is the coefficient of $\omega_0$ in the expansion of $\lambda$ into a 
linear combination of fundamental weights. But by definition, $-\alpha_0$ is 
the highest weight of ${\mathfrak{g}}_{-1}={\mathfrak{g}}_1^*$, and we obtain 
\begin{equation}\label{BBW1}
H^1({\mathfrak{g}}_{-1},{\mathbb{V}})=\Cart^k{\mathfrak{g}}_1\cart{\mathbb{E}}.
\end{equation}

There is a unique element in ${\mathfrak{g}}$ whose adjoint action is given by 
multiplication by $j$ on ${\mathfrak{g}}_j$ for $j=-1,0,1$, called the grading
element. The representation ${\mathbb{V}}$ splits into eigenspaces for the 
action of this element,
and it is convenient for our purposes to write
this decomposition as
$${\mathbb{V}}={\mathbb{V}}_0\oplus{\mathbb{V}}_1\oplus\cdots
\oplus{\mathbb{V}}_N,\quad\mbox{in which }{\mathbb{V}}_0={\mathbb{E}}
\mbox{ and }{\mathfrak{g}}_i{\mathbb{V}}_j\subseteq{\mathbb{V}}_{i+j}.$$
This is the algebraic source of (\ref{affinesplit}) and the number $N$ in
Theorems~\ref{specialaffineversion} and~\ref{generalversion}. The explicit
formulae for $N$ in \S\ref{formulation} may be obtained by observing that $N$
depends linearly on the coefficients of the fundamental weights in expressing
the highest weight and then verifying our formulae for the fundamental
representations. By construction, the homomorphisms $\partial$ and $\delta^*$
decrease and increase this grading on ${\mathbb{V}}$, respectively. 
Now (\ref{BBW0}) says that 
${\mathbb{V}}_i\stackrel{\partial}{\to}
{\mathfrak{g}}_1\otimes{\mathbb{V}}_{i-1}$ is injective~$\forall i\geq 1$.
The module $\Cart^k{\mathfrak{g}}_1\cart{\mathbb{E}}$ appears with multiplicity
one in ${\mathfrak{g}}_1\otimes{\mathbb{V}}$. Moreover, since
${\mathfrak{g}}_1$ increases the grading,
$\Cart^k{\mathfrak{g}}_1\cart{\mathbb{E}}$ resides in
${\mathfrak{g}}_1\otimes{\mathbb{V}}_{k-1}$. 
{From} (\ref{BBW1}), we conclude that
\begin{equation}\label{conclusion}
{\mathbb{V}}_i\stackrel{\partial}{\longhookrightarrow}
{\mathfrak{g}}_1\otimes{\mathbb{V}}_{i-1}
\stackrel{\partial}{\longrightarrow}
\Lambda^2{\mathfrak{g}}_1\otimes{\mathbb{V}}_{i-2}\quad
\mbox{is exact for }1\leq i\leq k-1\mbox{ and }i>k.\end{equation}
Now define $\phi_0:{\mathbb{V}}_0\to{\mathbb{E}}$ as the identity and
$\phi_i:{\mathbb{V}}_i\to\bigotimes^i{\mathfrak{g}}_1\otimes{\mathbb{E}}$
inductively as the composition:--
$$\textstyle{\mathbb{V}}_i\stackrel{\partial}{\longrightarrow}
{\mathfrak{g}}_1\otimes{\mathbb{V}}_{i-1}
\xrightarrow{{\rm{id}}\otimes\phi_{i-1}}
\bigotimes^i{\mathfrak{g}}_1\otimes{\mathbb{E}}.$$
Also set ${\mathbb{K}}=\ker:\bigodot^k{\mathfrak{g}}_1\otimes{\mathbb{E}}
\to\Cart^k{\mathfrak{g}}_1\cart{\mathbb{E}}$, the kernel of the Cartan
product.
\begin{lemma}\label{injectusw} The homomorphism
$\phi_i:{\mathbb{V}}_i\to\bigotimes^i{\mathfrak{g}}_1\otimes{\mathbb{E}}$
\newline {\rm (1)} is injective for all ~$i\geq 0$,
\newline {\rm (2)} has values in
$\bigodot^i{\mathfrak{g}}_1\otimes{\mathbb{E}}$,
\newline {\rm (3)} is an isomorphism
${\mathbb{V}}_i\stackrel{\simeq\quad}{\longrightarrow}
\bigodot^i{\mathfrak{g}}_1\otimes{\mathbb{E}}$, for $0\leq i\leq k-1$,
\newline {\rm (4)} is an isomorphism
${\mathbb{V}}_i\stackrel{\simeq\quad}{\longrightarrow}
(\bigodot^i{\mathfrak{g}}_1\otimes{\mathbb{E}})\cap
(\bigodot^{i-k}{\mathfrak{g}}_1\otimes{\mathbb{K}})$, for $i\geq k$.
\end{lemma}
\begin{proof} Statements (1)--(3) immediately follow by induction from
(\ref{conclusion}). When $i=k$, however, the sequence in (\ref{conclusion}) is
no longer exact. 
Rather, (\ref{BBW1}) implies that $\phi_k:
{\mathbb{V}}_k\hookrightarrow\bigodot^k{\mathfrak{g}}_1\otimes{\mathbb{E}}$ has
$\Cart^k{\mathfrak{g}}_1\cart{\mathbb{E}}$ as cokernel. This yields the
isomorphism
$\phi_k:{\mathbb{V}}_k\stackrel{\simeq\quad}{\longrightarrow}{\mathbb{K}}$,
which is (4) when $i=k$. For $i>k$ the exactness of (\ref{conclusion}) proves
(4) by induction.
\end{proof}
\noindent Let us denote by
$\phi_i^{-1}:\bigodot^i{\mathfrak{g}}_1\otimes{\mathbb{E}}\to{\mathbb{V}}_i$
the inverse of $\phi_i$ for $0\leq i\leq k-1$. Then, by construction and
since $\delta^*$ inverts $\partial$ on
$\im(\delta^*)={\mathbb{V}}_1\oplus\cdots\oplus{\mathbb{V}}_N$, we have:--
\begin{lemma}\label{coincides}
Although $\delta^*\circ({\rm{id}}\otimes\phi_{i-1}^{-1})$ is defined on
${\mathfrak{g}}_1\otimes\bigodot^{i-1}{\mathfrak{g}}_1\otimes{\mathbb{E}}$,
it coincides with $\phi_i^{-1}$ on
$\bigodot^i{\mathfrak{g}}_1\otimes{\mathbb{E}}$ for $1\leq i\leq k-1$.
\end{lemma}
\noindent We can also be more precise concerning the identification of
$1^{\mathrm{st}}$ cohomology in~(\ref{BBW1}). {From} the Hodge decomposition
(\ref{hodge}) and (\ref{congenial}), the endomorphism $\pi$ of
${\mathfrak{g}}_1\otimes{\mathbb{V}}$ given by
$\pi\ph=\ph-\delta^*\partial\ph- \partial\delta^*\ph$ is projection onto the
unique irreducible $G_0$-module isomorphic
to~$\Cart^k{\mathfrak{g}}_1\cart{\mathbb{E}}$. To fix this isomorphism, we take
\begin{equation}\label{fix}
\textstyle\Cart^k{\mathfrak{g}}_1\cart{\mathbb{E}}\longhookrightarrow
{\mathfrak{g}}_1\otimes\bigodot^{k-1}{\mathfrak{g}}_1\otimes{\mathbb{E}}
\xrightarrow{{\mathrm{id}}\otimes\phi_{k-1}^{-1}}
{\mathfrak{g}}_1\otimes{\mathbb{V}}_{k-1}
\stackrel{\pi}{\longrightarrow}\ker(\partial)\cap\ker(\delta^*).\end{equation}

\section{Proof of the main theorem}\label{proof}
To prove Theorem~\ref{generalversion}, we shall use the algebra of
\S\ref{algebra} as follows. Recall that $M$ is supposed to have a
$G_0'$-structure so any representation of $G_0'$ (and thus any representation
of $G_0$ or $G$ by restriction) induces an associated bundle on~$M$. Of course,
$E$ should be the bundle associated to ${\mathbb{E}}$ and we have already
observed that the bundle associated to ${\mathfrak{g}}_1$ is the bundle of
$1$-forms~$\Lambda^1$. Now we may transfer the constructions and conclusions
of~\S\ref{algebra} into geometry on~$M$. The $G$-module ${\mathbb{V}}$ induces
a graded vector bundle
$$\VB=\VB_0\oplus\VB_1\oplus\cdots\oplus\VB_N$$
on $M$. The complex (\ref{differential}) induces a complex of vector bundle
homomorphisms
\begin{equation}\label{complex}
0\longrightarrow\VB\stackrel{\partial}{\longrightarrow}
\Lambda^1\otimes\VB\stackrel{\partial}{\longrightarrow}
\Lambda^2\otimes\VB\stackrel{\partial}{\longrightarrow}\cdots\end{equation}
and, similarly, (\ref{codiff}) induces
\begin{equation}\label{cocomplex}
0\longrightarrow\VB\stackrel{\delta^*}{\longleftarrow}
\Lambda^1\otimes\VB\stackrel{\delta^*}{\longleftarrow}
\Lambda^2\otimes\VB\stackrel{\delta^*}{\longleftarrow}\cdots\end{equation}
so that
$E=\ker\partial:\VB\longrightarrow\Lambda^1\otimes\VB$ and (\ref{fix})
induces
\begin{equation}\label{sitting}\Cart^k\Lambda^1\cart E\cong
\frac{\ker\partial:\Lambda^1\otimes\VB\longrightarrow
\Lambda^2\otimes\VB}
{\mbox{$\im\partial:\VB\longrightarrow\Lambda^1\otimes\VB$}}=
\ker(\partial)\cap\ker(\delta^*).\end{equation}
Lemma~\ref{injectusw} part (3) yields
$$\textstyle\phi_j:\VB_i\stackrel{\simeq\quad}{\longrightarrow}
\bigodot^j\Lambda^1\otimes E\quad\mbox{for}\quad 0\leq j\leq k-1.$$
Lemma~\ref{injectusw} part (4) identifies $V_i$ with the classical
prolongations (\ref{classicalV}) for $i\geq k$.

\subsection*{A splitting operator}
According to the statement of Theorem~\ref{generalversion} we should choose a
connection $\nabla$ on $M$ that is compatible with the $G_0'$-structure. {From}
this, we obtain connections on all associated vector bundles, in particular on
$E$ and $\VB$. Being induced from a principal $G_0'$-connection, they respect
the grading on $\VB$ and commute with the homomorphisms in (\ref{complex})
and~(\ref{cocomplex}). We shall denote all of these linear connections
by~$\nabla$.

To prove Theorem~\ref{generalversion} we shall construct $L:E=\VB_0\to\VB$,
an $N^{\mathrm{th}}$ order linear differential operator, so that
$\sigma\mapsto L\sigma$ induces the isomorphism (\ref{genisomorphism}). Since
the isomorphism in the other direction should simply be given by
$\sigma=\Sigma_0$, the component of $\Sigma$ in $\VB_0=E$, the composition
$\sigma\mapsto (L\sigma)_0$ should be the identity. For this reason we refer
to $L$ as a `splitting operator'. Its definition is
\begin{equation}\label{defL}
L\sigma=\sum_{i=0}^N(-1)^i(\delta^*\circ\nabla)^i\sigma.\end{equation}
Of course, this an $N^{\mathrm{th}}$~order linear differential operator.
Moreover, since $\sigma$ is a section of $E=\VB_0$, we see that
$(\delta^*\circ\nabla)^i\sigma$ is a section of $\VB_i$ and that a section
$\Sigma=(\Sigma_0,\Sigma_1,\dots,\Sigma_N)$ of $\VB$ is of the form $L\sigma$
if and only if
\begin{equation}\label{rangeL}\Sigma_0=\sigma\quad\mbox{and}\quad
\Sigma_i=-\delta^*\nabla\Sigma_{i-1}\mbox{ for }1\leq i\leq N.\end{equation}

Next, we define the connection $\widetilde\nabla$ on $\VB$ as
$\widetilde\nabla=\nabla+\partial$. So we simply add the algebraic operator
$\partial:\VB\to\Lambda^1\otimes\VB$ to the component-wise connection
$\nabla$. Of course, this defines a linear connection. Note, however, that
whilst for a section $\Sigma_i$ of $\VB_i$, the covariant derivative
$\nabla\Sigma_i$ is a $1$-form with coefficients in~$\VB_i$, the algebraic term
$\partial\Sigma_i$ is a $1$-form with coefficients in~$\VB_{i-1}$. Otherwise
put, for a section $\Sigma=(\Sigma_0,\Sigma_1,\dots,\Sigma_N)$ of $\VB$ the
component of $\widetilde\nabla\Sigma$ taking values in $\Lambda^1\otimes\VB_i$
is $\nabla\Sigma_i+\partial\Sigma_{i+1}$ for $i<N$, while for $i=N$ we simply
obtain $\nabla\Sigma_N$.

We should now compute the curvature of~$\widetilde\nabla$. The curvatures of
all the connections $\nabla$ are induced by the same $2$-form~$R$, which acts
on the sections of any associated bundle. On the other hand, for the induced
connection on $TM$ we also have the torsion, which we view as a section of
$\Lambda^2\otimes TM$. (In the affine and Riemannian cases we can always choose
$\nabla$ to be torsion-free but not with a general $G_0'$-structure).
\begin{lemma}\label{curvaturelemma}
Let $R$ be the curvature of the connections $\nabla$ and $T$ the
torsion of the connection $\nabla$ on~$TM$. Let
$\widetilde R\in\Gamma(\Lambda^2\otimes\End(\VB,\VB))$ be the curvature
of~$\widetilde\nabla$. Then for vector fields $\xi$ and $\eta$ on $M$ and a
section $\Sigma=(\Sigma_0,\dots,\Sigma_N)$ of~$\VB$, the $\VB_i$-component
of $\widetilde R(\xi,\eta)\Sigma$ is given by
$$R(\xi,\eta)\Sigma_i+(\partial\Sigma_{i+1})(T(\xi,\eta)).$$
In particular, $\widetilde\nabla$ is flat if and only if $\nabla$ has zero
curvature and torsion.
\end{lemma}
\begin{proof}
By definition, $\widetilde\nabla_\xi\widetilde\nabla_\eta\Sigma=
\widetilde\nabla_\xi(\nabla_\eta\Sigma+(\partial\Sigma)(\eta))$. Writing out
the first operator as $\nabla+\partial$, we obtain
\begin{equation}\label{comp}
\nabla_\xi\nabla_\eta\Sigma+\nabla_\xi((\partial\Sigma)(\eta))+
(\partial(\nabla_\eta\Sigma))(\xi)+(\partial(\partial\Sigma)(\eta))(\xi).
\end{equation}
To obtain $\widetilde R(\xi,\eta)\Sigma$ we should subtract the same sum
with $\xi$ and $\eta$ exchanged and then subtract
\begin{equation}\label{tor}
\widetilde\nabla_{[\xi,\eta]}\Sigma
=\nabla_{[\xi,\eta]}\Sigma+(\partial\Sigma)([\xi,\eta]).
\end{equation}
On the Lie algebra level $(\partial v)(Y)=Yv$ and thus
$\partial((\partial v)(Y))(Z)=Z(Yv)$, which is symmetric in $Y$ and $Z$ since
${\mathfrak{g}}_1$ is an Abelian Lie algebra. Hence the last term in
(\ref{comp}) vanishes after exchange and subtraction. Also, we may write
$$\nabla_\xi((\partial\Sigma)(\eta))=
(\nabla_\xi(\partial\Sigma))(\eta)+(\partial\Sigma)(\nabla_{\xi}\eta)$$
and, since $\partial$ is parallel, rewrite the first summand as
$(\partial(\nabla_\xi\Sigma))(\eta)$. But this cancels with one of the terms
from the other summand of the form (\ref{comp}). Altogether, we see that the
last three terms in the two summands of the form (\ref{comp}) together
contribute $(\partial\Sigma)(\nabla_{\xi}\eta-\nabla_\eta\xi)$. Subtracting the
last term in (\ref{tor}) we obtain $(\partial\Sigma)(T(\xi,\eta))$ by
definition of the torsion. On the other hand, the first terms in the two
summands of the form (\ref{comp}) add up with the remaining term of (\ref{tor})
to $R(\xi,\eta)\Sigma$. Now, the result follows by splitting into components.
\end{proof}

Having at hand the operators $L$ and $\widetilde\nabla$, we now define an
operator $E=\VB_0\to\Lambda^1\otimes\VB$ as the composition
$\widetilde\nabla\circ L$. {From} (\ref{sitting}) we know that
$F=\Cart^k\Lambda^1\cart E$ sits as the subbundle
$\ker(\partial)\cap\ker(\delta^*)$ in $\Lambda^1\otimes\VB$, and we can use the
algebraic Hodge structure to define a projection onto this subbundle. Indeed,
in \S\ref{algebra} we arranged that this projection be explicitly given by
$\ph\mapsto \pi\ph\equiv\ph-\delta^*\partial\ph-
\partial\delta^*\ph$.
Using~(\ref{fix}), we now define a differential operator $D^\nabla:E\to F$ by
$D^\nabla\equiv
(-1)^{k-1}({\rm{id}}\otimes\phi_{k-1})\circ\pi\circ\widetilde\nabla\circ L$.
The main properties of $L$ and $D^\nabla$ are collected in:--

\begin{proposition}\label{properties}\mbox{ }
\newline {\rm (1)} A section $\Sigma=(\Sigma_0,\dots,\Sigma_N)$ of $\VB$ lies
in the image of $L$ if and only if $\delta^*(\widetilde\nabla\Sigma)=0$ and, if
this is the case, then $\Sigma=L(\Sigma_0)$.
\newline {\rm (2)} Mapping $\sigma\in\Gamma(E)$ to the components of
$L\sigma$ in $\VB_0\oplus\dots\oplus\VB_i$ induces a vector bundle
homomorphism $J^i\VB_0\to\VB_0\oplus\dots\oplus\VB_i$, which is an isomorphism
for~$i<k$.
\newline {\rm (3)} The differential operator $D^\nabla:E\to F$ is of order $k$
and its symbol is the Cartan product.
\end{proposition}
\begin{proof}
(1)  Since $\VB\stackrel{\delta^*}{\longleftarrow}\Lambda^1\otimes\VB$ inverts
$\partial$ on~$\im(\partial)$, we may easily compute the components
of~$\delta^*(\widetilde\nabla\Sigma)$. We find that
$\delta^*(\widetilde\nabla\Sigma)_0=0$ and, for $1\leq i\leq N$,
$$\delta^*(\widetilde\nabla\Sigma)_i=\delta^*((\widetilde\nabla\Sigma)_{i-1})
=\delta^*(\nabla\Sigma_{i-1}+\partial\Sigma_i)=
\delta^*(\nabla\Sigma_{i-1})+\Sigma_i$$
whose vanishing is exactly the criterion (\ref{rangeL}) we already found for
$\Sigma=(\Sigma_0,\dots\Sigma_N)$ to be in the range of~$L$. In (\ref{rangeL})
we also observed that, in this case, $\Sigma=L(\Sigma_0)$.

\noindent (2) By construction, mapping $\sigma$ to the $\VB_i$-component of
$L\sigma$ is a linear differential operator of order at most~$i$. Thus, we
obtain $J^i\VB_0\to\VB_0\oplus\dots\oplus\VB_i$ for all $i=0,\dots,N$. We can
compute the leading terms of $(L\sigma)_i$ quite explicitly as follows.
Firstly, $(L\sigma)_0$ is just~$\sigma$, a section of $\VB_0=E$. Next, from its
definition~(\ref{defL}), we have $(L\sigma)_1=-\delta^*\nabla\sigma$. Assuming
that $1<k$, we see from Lemma~\ref{coincides} that $\delta^*$ coincides
with~$\phi_1^{-1}$. Therefore, $(L\sigma)_1=-\phi_1^{-1}\nabla\sigma$, a
section of $\VB_1$. Now
$\nabla(L\sigma)_1=-({\rm{id}}\otimes\phi_1^{-1})\nabla(\nabla\sigma)$,
where $\nabla(\nabla\sigma)$ is a section of
$\Lambda^1\otimes\Lambda^1\otimes E$. But if we decompose
$$\textstyle\Lambda^1\otimes\Lambda^1\otimes E=
(\bigodot^2\Lambda^1\otimes E)\oplus(\Lambda^2\otimes E),$$
then the component $\nabla\wedge\nabla\sigma$ of $\nabla(\nabla\sigma)$ is a
zeroth order operator (made from curvature and torsion). If $2<k$, then from
Lemma~\ref{coincides} we conclude that
$$(L\sigma)_2=-\delta^*\nabla(L\sigma)_1
=\delta^*({\rm{id}}\otimes\phi_1^{-1})\nabla(\nabla\sigma)=
\phi_2^{-1}\nabla\odot\nabla\sigma+\mbox{lots},$$
where `${\mathrm{lots}}$' stands for `lower order terms' (in this case zeroth
order). By induction, we claim that
\begin{equation}\label{claim}(L\sigma)_i=
(-1)^i\phi_i^{-1}\underbrace{\nabla\odot\nabla\odot\cdots\odot\nabla}_i\sigma
+\mbox{lots},\quad\mbox{for }0\leq i \leq k-1.\end{equation}
For the inductive step, observe that
$$\nabla_a\nabla_{(b}\nabla_c\cdots\nabla_{d)}=
\nabla_{(a}\nabla_b\nabla_c\cdots\nabla_{d)}+{\mathrm{lots}}$$
as differential operators. Therefore,
$$\begin{array}{rcl}\nabla(L\sigma)_{i-1}&=&\nabla((-1)^{i-1}\phi_{i-1}^{-1}
\underbrace{\nabla\odot\nabla\odot\cdots\odot\nabla}_{i-1}\sigma
+\mbox{lots})\\[17pt]
&=&(-1)^{i-1}({\rm{id}}\otimes\phi_{i-1}^{-1})
(\underbrace{\nabla\odot\nabla\odot\nabla\odot\cdots\odot\nabla}_i\sigma
+\mbox{lots})\end{array}$$
and so, for $i<k$,
$$\begin{array}{rcccl}(L\sigma)_i&=&-\delta^*\nabla(L\sigma)_{i-1}
&=&(-1)^i\delta^*({\rm{id}}\otimes\phi_{i-1}^{-1})
(\underbrace{\nabla\odot\nabla\odot\nabla\odot\cdots\odot\nabla}_i\sigma
+\mbox{lots})\\[17pt]
&&&=&(-1)^i\phi_i^{-1}
\underbrace{\nabla\odot\nabla\odot\nabla\odot\cdots\odot\nabla}_i\sigma
+\mbox{lots},
\end{array}$$
the last equality coming from Lemma~\ref{coincides}. We have shown
(\ref{claim}) and, clearly, this is sufficient to establish~(2).

\noindent (3) The projection
$$\Lambda^1\otimes E\ni
\ph\mapsto\pi\ph\equiv\ph-\delta^*\partial\ph-\partial\delta^*\ph
\in\ker(\partial)\cap\ker(\delta^*)$$
kills $\im(\partial)$ so $D^\nabla\sigma=
(-1)^{k-1}({\rm{id}}\otimes\phi_{k-1})(\pi(\nabla(L\sigma)_{k-1}))$.
{From} (\ref{fix}) and (\ref{claim}) we see that
$$D^\nabla\sigma=
\pi(\nabla(\underbrace{\nabla\odot\nabla\odot\cdots\odot\nabla}_{k-1}\sigma
+{\mathrm{lots}})),$$
where now $\pi:\Lambda^1\otimes\bigodot^{k-1}\Lambda^1\otimes E
\to\Cart^k\Lambda^1\otimes E =F$ denotes canonical projection onto this
irreducible tensor bundle. It is now clear the $D^\nabla$ has the Cartan
product as its symbol. \end{proof}

\subsection*{First step}
Now we can perform the first step in rewriting the equation $D\sigma=0$
on sections of $E$ in terms of sections of $\VB$:--
\begin{proposition}
Let $D:E\to F$ be a $k^{\mathrm{th}}$ order semilinear differential operator
as in Theorem~\ref{generalversion}. Then there is a fibre bundle homomorphism
$A:\VB_0\oplus\dots\oplus\VB_{k-1}\to F$ such that $\sigma\mapsto L\sigma$
induces a set bijection
$$\{\sigma\in\Gamma(E)\mbox{\rm\ s.t. }D\sigma=0\}\cong
\{\Sigma\in\Gamma(\VB)\mbox{\rm\ s.t. }
\widetilde\nabla\Sigma+A(\Sigma)\in\Gamma(\im(\delta^*))\}.$$ If
$D$ is linear, then $A$ is linear, i.e.~a vector bundle homomorphism.
\end{proposition}
\begin{proof}
{From} part (3) of Proposition~\ref{properties} we conclude that the operators
$D$ and $D^\nabla$ have the same symbol. Therefore, we may write
$D\sigma=D^\nabla\sigma+\Psi(j^{k-1}\sigma)$ for some bundle map
$\Psi:J^{k-1}E\to F$. By part (2) of Proposition~\ref{properties} there is a
unique fibre bundle map $A:\VB_0\oplus\dots\oplus\VB_{k-1}\to F$ (which we may
extend trivially to~$\VB$) such that $\Psi(j^{k-1}\sigma)=(-1)^{k-1}A(L\sigma)$
for all $\sigma\in\Gamma(E)$. Of course, if $D$ is linear, then $\Psi$ is a
vector bundle homomorphism and hence $A$ is a vector bundle homomorphism too.

Now $\widetilde\nabla L\sigma$ is a section of $\ker(\delta^*)$ by part (1)
of Proposition~\ref{properties} and the same is true for $A(L\sigma)$ since, by
construction, $A$ even has values in~$F=\ker(\partial)\cap\ker(\delta^*)$.
The last observation even shows that $\pi(A(L\sigma))=A(L\sigma)$ for
any~$\sigma$. Hence, vanishing of $D\sigma=(-1)^{k-1}\pi(\widetilde\nabla
L\sigma+A(L\sigma))$ is equivalent to $\widetilde\nabla L\sigma+A(L\sigma)$
being a section of the subbundle $\im(\delta^*)$.

Conversely, assume that $\Sigma\in\Gamma(\VB)$ has the property that
$\widetilde\nabla\Sigma+A(\Sigma)$ is a section of $\im(\delta^*)$. Then, in
particular it is a section of $\ker(\delta^*)$ and since $\delta^*(A(\Sigma))$
always vanishes we conclude that $\delta^*(\widetilde\nabla\Sigma)=0$. By
part (1) of Proposition~\ref{properties} this implies $\Sigma=L(\Sigma_0)$
and, as above, we see that $D(\Sigma_0)=0$.
\end{proof}

\subsection*{Second step}
The next step in the procedure is to show that, if
$\widetilde\nabla\Sigma+A(\Sigma)$ is a section of $\im(\delta^*)$, then its
value can be actually computed. We shall do this in a more general situation
than needed for the proof of Theorem~\ref{generalversion}. The motivation for
this is that if $A$ is linear, then it can be absorbed into the connection, so
dealing with a more general class of connections is helpful. Notice that any
smooth section of the bundle $\im(\delta^*)\subset\Lambda^1\otimes\VB$ can be
written as $\delta^*\psi$ for some smooth $\psi\in\Gamma(\Lambda^2\otimes\VB)$.

\begin{proposition}\label{1.5}
Let $\nabla$ be a linear connection on $\VB$ such that for each $i=0,\dots,N$
and each smooth section $\Sigma\in\Gamma(\VB)$ that has values in $\VB_i$ only,
the covariant derivative $\nabla\Sigma$ lies in
$\Gamma(\Lambda^1\otimes(\VB_i\oplus\dots\oplus\VB_N))$ and put
$\widetilde\nabla=\nabla+\partial$. Let $A:\VB\to\Lambda^1\otimes\VB$ be a
fibre bundle map such that for $v=(v_0,v_1,\ldots,v_N)\in V$ the component of
$A(v)$ in $\Lambda^1\otimes V_i$ depends only on $v_0,\ldots,v_i$. Then there
is a fibre bundle map
$$B:J^N\VB=J^N\VB_0\oplus\dots\oplus J^N\VB_N\to
\Lambda^1\otimes\VB$$
such that $\widetilde\nabla\Sigma+A(\Sigma)\in\Gamma(\im(\delta^*))$ is
equivalent to $\widetilde\nabla\Sigma+B(j^N\Sigma)=0$. Moreover, the component
$B_i$ of $B$ with values in $\Lambda^1\otimes\VB_i$ factors through
$$J^i\VB_0\oplus J^{i-1}\VB_1\oplus\dots\oplus J^1\VB_{i-1}\oplus\VB_i.$$
If $A$ is linear then $B$ can be chosen to be a vector bundle homomorphism.
\end{proposition}
\begin{proof}
Suppose that $\widetilde\nabla\Sigma+A(\Sigma)+\delta^*\psi=0$ for some
$\psi\in\Gamma(\Lambda^2\otimes\VB)$. Recall that the linear connection
$\widetilde\nabla$ on $\VB$ extends to an operation $d^{\widetilde\nabla}$
on $\VB$-valued forms called the covariant exterior derivative. For
$\alpha\in\Gamma(\Lambda^1\otimes\VB)$ the covariant exterior derivative is
explicitly given by
$$ d^{\widetilde\nabla}\alpha(\xi,\eta)=\widetilde\nabla_\xi(\alpha(\eta))-
\widetilde\nabla_\eta(\alpha(\xi))-\alpha([\xi,\eta]),$$
for all vector fields $\xi$ and $\eta$ on $M$. Clearly, $d^{\widetilde\nabla}$
is a first order differential operator. Moreover, if
$\alpha=\widetilde\nabla\Sigma$ for some $\Sigma\in\Gamma(\VB)$, then this
definition immediately implies that
$d^{\widetilde\nabla}\widetilde\nabla\Sigma(\xi,\eta)=
\widetilde R(\xi,\eta)(\Sigma)$.

Now we define $B$ inductively as follows. We put
$B_0(\Sigma)\equiv A_0(\Sigma)$. By assumption, this is algebraic (i.e.\ of
order zero) in $\Sigma$ and depends only on the component $\Sigma_0$. Let
$\widetilde R\bullet\Sigma$ denote the $\VB$-valued $2$-form
$(\xi,\eta)\mapsto \widetilde R(\xi,\eta)(\Sigma)$. Having defined the
components $B_j$ for $j<i$, take the component
$(\widetilde R\bullet\Sigma+
d^{\widetilde\nabla}(B_{i-1}(\Sigma)+\dots+B_0(\Sigma)))_{i-1}$ in
$\Gamma(\Lambda^2\otimes\VB_{i-1})$ and define
\begin{equation}\label{defB}B_i(\Sigma)\equiv A_i(\Sigma)
-\delta^*\left((\widetilde R\bullet\Sigma
+d^{\widetilde\nabla}(B_{i-1}(\Sigma)+\dots+B_0(\Sigma)))_{i-1}
+\partial(A_i(\Sigma))\right).\end{equation}
By assumption, $A$ is algebraic in $\Sigma$ and $A_i(\Sigma)$ depends only on
the components $\Sigma_0,\dots,\Sigma_i$. To understand the dependence of
$\widetilde R\bullet\Sigma$, note that by assumption on $\nabla$, the form
$(\widetilde\nabla\Sigma)_j$ depends only on $\Sigma_0,\dots,\Sigma_{j+1}$.
Hence the $\VB_j$-component of $\widetilde R(\xi,\eta)(\Sigma)$ depends at most
on $\Sigma_0,\dots,\Sigma_{j+2}$ (since computing curvature needs two
derivatives). However, as in the proof of Lemma \ref{curvaturelemma}, we see
that for $\Sigma\in\Gamma(\VB_{j+2})$ the only contribution of
$\widetilde R(\xi,\eta)(\Sigma)$ in $\VB_j$ is
$\partial((\partial\Sigma)(\eta))(\xi)-\partial((\partial\Sigma)(\xi))(\eta)$
and we have shown that this vanishes. Hence, the term
$(\widetilde R\bullet\Sigma)_{i-1}$ depends only on $\Sigma_0,\dots,\Sigma_i$.
Assuming inductively that for $\ell\leq i-1$, the value $B_\ell(\Sigma)(x)$
depends only on $j_x^\ell\Sigma_0,j_x^{\ell-1}\Sigma_1,\dots,\Sigma_\ell(x)$
for each $x\in M$, we immediately conclude from the fact that
$d^{\widetilde\nabla}$ is first order that $B_i(\Sigma)(x)$ depends only on
$j_x^i\Sigma_0, j_x^{i-1}\Sigma_1,\dots,\Sigma_i(x)$. Hence
our components $B_i$ define a bundle map $B$ whose dependence on jets is
exactly as required. Moreover, if $A$ is linear, then obviously $B$ is a
vector bundle homomorphism.

Next we show that the equation $\widetilde\nabla\Sigma+B(\Sigma)=0$ is
equivalent to $\widetilde\nabla\Sigma+A(\Sigma)$ being a section
of~$\im(\delta^*)$. On the one hand, we see from the definition (\ref{defB})
that $A(\Sigma)-B(\Sigma)$ is a section of $\im(\delta^*)$ for any
$\Sigma\in\Gamma(\VB)$. Thus $\widetilde\nabla\Sigma+B(\Sigma)=0$ implies that
$\widetilde\nabla\Sigma+A(\Sigma)$ has values in $\im(\delta^*)$.

Conversely, assume that $\widetilde\nabla\Sigma+A(\Sigma)+\delta^*\psi=0$ for
some $\psi\in\Gamma(\Lambda^2\otimes\VB)$. Then we claim that
$A(\Sigma)+\delta^*\psi=B(\Sigma)$. Since $\delta^*$ has values in
$\Lambda^1\otimes(\VB_1\oplus\dots\oplus\VB_N)$ and, by definition,
$B_0(\Sigma)=A_0(\Sigma)$, this is true for the component in
$\Gamma(\Lambda^1\otimes\VB_0)$.

To proceed inductively, we need one more observation concerning
$d^{\widetilde\nabla}$. Suppose that
$\alpha\in\Gamma(\Lambda^1\otimes(\VB_i\oplus\dots\oplus\VB_N))$. Then, from
the formula for $d^{\widetilde\nabla}$, it is manifest that
$d^{\widetilde\nabla}\alpha\in
\Gamma(\Lambda^2\otimes(\VB_{i-1}\oplus\dots\oplus\VB_N)$ and the component
$(d^{\widetilde\nabla}\alpha)_{i-1}$ is easy to compute: expanding
$\widetilde\nabla=\nabla+\partial$ in the above formula, we see that
$$(d^{\widetilde\nabla}\alpha)_{i-1}(\xi,\eta)=\partial(\alpha_i(\eta))(\xi)-
\partial(\alpha_i(\xi))(\eta)$$
and, looking at the definition of $\partial$, this means that
$(d^{\widetilde\nabla}\alpha)_{i-1}=\partial(\alpha_i)$.

Now suppose inductively that $(A(\Sigma)+\delta^*\psi)_\ell=B_\ell(\Sigma)$
for $\ell=0,\dots,i-1$. Denoting by the subscript $\geq i$ the components with
values in $\Lambda^1\otimes(\VB_i\oplus\dots\oplus\VB_N)$ we may rewrite the
equation $\widetilde\nabla\Sigma+A(\Sigma)+\delta^*\psi=0$ as
$$\widetilde\nabla\Sigma+B_0(\Sigma)+\dots+B_{i-1}(\Sigma)+A_{\geq i}(\Sigma)+
(\delta^*\psi)_{\geq i}=0.$$
Applying $d^{\widetilde\nabla}$ and looking at the component in
$\Lambda^2\otimes\VB_{i-1}$ we obtain
$$0=(\widetilde R\bullet\Sigma+
d^{\widetilde\nabla}(B_0(\Sigma)+\dots+B_{i-1}(\Sigma)))_{i-1}+
\partial(A_i(\Sigma))+\partial(\delta^*\psi)_i.$$
Applying $\delta^*$, the last term gives $(\delta^*\psi)_i$
and from (\ref{defB}) we see $A_i(\Sigma)+(\delta^*\psi)_i=B_i(\Sigma)$, which
completes the proof.
\end{proof}

\subsection*{Third step}
The final reduction is now done by solving component by component:--
\begin{proposition}
Suppose that $\nabla$ is a connection on $\VB$ satisfying the hypothesis of
Proposition \ref{1.5} and
$$B:J^N\VB=J^N\VB_0\oplus\dots\oplus J^N\VB_N\to
\Lambda^1\otimes\VB$$
is a fibre bundle map such that the component $B_i$ of $B$ in
$T^*M\otimes\VB_i$ factors through
$J^i\VB_0\oplus J^{i-1}\VB_1\oplus\dots\oplus J^1\VB_{i-1}\oplus\VB_i$.

Then there is a fibre bundle map $C:\VB\to\Lambda^1\otimes\VB$ such
that $\widetilde\nabla\Sigma+B(\Sigma)=0$ is equivalent to
$\widetilde\nabla\Sigma+C(\Sigma)=0$. If $B$ is a vector bundle homomorphism,
then also $C$ can be chosen to be a vector bundle homomorphism.
\end{proposition}
\begin{proof}
Choosing a connection on~$TM$, we may form iterated covariant derivatives of
sections of $\VB$ and by the assumptions on $B$ we may write the components of
$B$ (with the obvious meaning of subscripts) as
$$B_i(\Sigma)=B_i(\Sigma_{\leq i},(\widetilde\nabla\Sigma)_{\leq i-1},
\dots,(\widetilde\nabla^i\Sigma)_0).$$
The component in $\Lambda^1\otimes\VB_0$ of $\widetilde\nabla\Sigma+B(\Sigma)$
is given by $(\widetilde\nabla\Sigma)_0+B_0(\Sigma_0)$ and we simply put
$C_0(\Sigma)\equiv B_0(\Sigma_0)$. The next component has the form
$(\widetilde\nabla\Sigma)_1+B_1(\Sigma_0,\Sigma_1,(\widetilde\nabla\Sigma)_0)$.
Defining $C_1(\Sigma_0,\Sigma_1)\equiv B_1(\Sigma_0,\Sigma_1,-C_0(\Sigma_0))$,
we see that vanishing of $(\widetilde\nabla\Sigma+B(\Sigma))_{\leq 1}$ is
equivalent to vanishing of $(\widetilde\nabla\Sigma+C(\Sigma))_{\leq 1}$, where
$C=C_0+C_1$.

Let us inductively assume that $i>1$ and we have found a fibre bundle map
$C:\VB\to\Lambda^1\otimes(\VB_0\oplus\dots\oplus\VB_{i-1})$ such that vanishing
of $(\widetilde\nabla\Sigma+B(\Sigma))_{\leq i-1}$ is equivalent to vanishing
of $(\widetilde\nabla\Sigma+C(\Sigma))_{\leq i-1}$ and such that the component
$C_j(\Sigma)$ depends only on $\Sigma_0,\dots,\Sigma_j$ for each $j<i$. Let us
also assume that we have derived, for any $\Sigma$ such that
$(\widetilde\nabla\Sigma+C(\Sigma))_{\leq i-1}=0$, formulae for
$(\widetilde\nabla^\ell\Sigma)_{\leq i-\ell}$ as algebraic expressions
in $\Sigma_{\leq i}$.

So by assumption we have formulae for all the terms going into $B_i$ as
algebraic operators in $\Sigma_{\leq i}$ and inserting these formulae, we
obtain a bundle map $C_i$ with values in $\Lambda^1\otimes\VB_i$, which depends
only on $\Sigma_{\leq i}$. By construction, vanishing of
$(\widetilde\nabla\Sigma+B(\Sigma))_{\leq i}$ is equivalent to vanishing of
$(\widetilde\nabla\Sigma+C(\Sigma))_{\leq i}$. Suppose now that $\Sigma$
satisfies this equation. By the assumption on $\widetilde\nabla$, vanishing of
$(\widetilde\nabla\Sigma+C(\Sigma))_{\leq i}$ implies vanishing of
$(\widetilde\nabla^{\ell}(\widetilde\nabla\Sigma+C(\Sigma)))_{\leq i-\ell}$
for each $\ell=1,\dots,i$. Similarly,
$(\widetilde\nabla^{\ell}(C(\Sigma)))_{\leq i-\ell}$ depends algebraically on
$C(\Sigma)_{\leq i}$, to first order on $C(\Sigma)_{\leq i-1}$ and so on.
Hence, expanding this, it can be written as an expression in $\Sigma_{\leq i}$,
$(\widetilde\nabla\Sigma)_{\leq i-1}$,\dots,
$(\widetilde\nabla^\ell\Sigma)_{\leq i-\ell}$ and we have algebraic formulae
for all these by inductive hypothesis. Thus we see that vanishing of
$(\widetilde\nabla^{\ell}(\widetilde\nabla\Sigma+C(\Sigma)))_{\leq i-\ell}$
gives us an algebraic expression for
$(\widetilde\nabla^{\ell+1}\Sigma)_{\leq i-\ell}$ for each $\ell=1,\dots,i$,
which completes the inductive step. Of course, linearity is never lost in this
process, so if one starts with a linear operator~$B$, one will end up with a
vector bundle homomorphism~$C$. \end{proof}

Since the output of each step of our rewriting procedure is a special
case of the input of the next step, this completes the proof of
Theorem~\ref{generalversion}.

\begin{remark} As far as the proof of Theorem~\ref{generalversion} is
concerned, the only r\^ole that $\partial^*$ played was in constructing
$\delta^*$ as a left inverse to $\partial$. Of course, the definition of
$\partial^*$ and the resulting algebraic Hodge theory is extremely natural but,
in defining $\partial$, only the structure of ${\mathbb{V}}$ as a
${\mathfrak{g}}_{-1}$-module is needed. It is also important that $\partial$
respect the $G_0$-action but, as far as $\delta^*$ goes, any other
$G_0$-invariant splittings would work just as well. In practise, there can be
considerably simpler {\em ad hoc\/} choices.
\end{remark}

\begin{remark} The dimension bound of Corollary~\ref{dimbound} is sharp. The
bound is attained by choosing a manifold $M$ endowed with a $G_0'$-structure
and a compatible connection, such that all the connections $\nabla$ have zero
curvature and the connection $\nabla$ on $TM$ also has zero torsion. Such an
example is always provided by the constant $G_0'$-structure on ${\mathbb{R}}^n$
(where $n=\dim({\mathfrak{g}}_{-1})$) with the standard flat connection. In
this case, let us consider the equation $D^\nabla\sigma=0$. Then our first step
of rewriting simply leads to $\widetilde\nabla\Sigma+\delta^*\psi=0$.
Applying $\delta^*d^{\widetilde\nabla}$, the first term does not give any
contribution, since $\widetilde\nabla$ has zero curvature by
Lemma~\ref{curvaturelemma}. This implies that $\delta^*\psi=0$. Hence the
whole rewriting is already finished and we conclude that the differential
splitting $L:E\to\VB$ induces a bijection between solutions of
$D^{\nabla}\sigma=0$ and sections $\Sigma\in\Gamma(\VB)$ that are parallel for
the flat connection $\widetilde\nabla$. Locally, a flat connection always
has the maximal dimension for its space of parallel sections.
\end{remark}

\begin{remark}
The flat case also shows that the bound $N$ on the order of the jet of $\sigma$
at a point $p\in M$ needed uniquely to specify a solution of $D^\nabla\sigma=0$
is sharp. To see this, note that $\widetilde\nabla\Sigma=0$ in the flat case is
equivalent to $\nabla\Sigma_i=-\partial\Sigma_{i+1}$, for all~$i$. Therefore,
$$\Sigma|_p\in (V_N)_p\Rightarrow\nabla\Sigma|_p\in (V_{N-1}\oplus V_N)_p
\Rightarrow\cdots\Rightarrow\nabla^{N-1}\Sigma|_p\in
(V_1\oplus\cdots\oplus V_N)_p,$$
whence $\nabla^{N-1}\sigma|_p=\nabla^{N-1}\Sigma_0|_p=0$. But, since
$\widetilde\nabla$ is flat, there is no problem finding a parallel section
$\Sigma$ of $V$ with $\Sigma|_p$ lying in~$(V_N)_p$.
\end{remark}

\begin{remark}
In this flat case, the operator $D^\nabla$ is the first in the so-called
`Bernstein-Gelfand-Gelfand (BGG) resolution' and one motivation for our study
comes from analogues of these first operators on almost Hermitian symmetric
manifolds~\cite{baston} or, more generally, on parabolic
geometries~\cite{cd,cssannals}. By construction, these analogues are invariant
linear differential operators having the same symbol as in the flat case. The
$G_0'$-geometries studied in this article cover the almost Hermitian symmetric
case so Theorem~\ref{generalversion} covers the first BGG operators on
these geometries. This includes the various so-called `conformal Killing' or
`twistor' equations in conformal geometry.
\end{remark}

\begin{remark} A useful viewpoint on the outcome of
Theorem~\ref{generalversion} is that it restricts the possible jets of $\sigma$
that might be specified at a point for a solution of~$D\sigma=0$. In the flat
case and the equation $D^\nabla\sigma=0$, these jets may be freely specified.
In general, there are further constraints, which may be obtained by
cross-differentiation of the closed system
$\widetilde\nabla\Sigma+\Phi(\Sigma)=0$.
\end{remark}

\end{document}